\documentclass[11pt]{article}
\usepackage{graphicx}
\usepackage{amsmath,amssymb,amsthm,amsfonts,mathrsfs}
\usepackage{amssymb}
\usepackage{enumerate}
\usepackage[numbers,sort&compress]{natbib}

\newtheorem{thm}{Theorem}[section]

\newtheorem{lem}[thm]{Lemma}

\theoremstyle{definition}
\newtheorem{defn}{Definition}[section]
\theoremstyle{remark}
\newtheorem{Rem}{Remark}[section]

\numberwithin{equation}{section}

\DeclareMathSymbol{\C}{\mathalpha}{AMSb}{"43}

\newcommand{\N}{{\mathbb N}}

\newcommand{\R}{{\mathbb{R}}}

\newcommand{\bsub}{\begin{subequations}}
\newcommand{\esub}{\end{subequations}$\!$}

\begin{document}

\title{Normalized solutions for logarithmic Schr\"{o}dinger equation with a perturbation of power law nonlinearity\thanks {The research was supported by the National Natural Science Foundation of China (12071170).
}
}

\author{ Wei Shuai\thanks{School of Mathematics and Statistics, \& Hubei Key Laboratory of Mathematical Sciences, Central China Normal University, Wuhan 430079, P. R. China.
Email: \texttt{wshuai@mail.ccnu.edu.cn}.}
\ \ \ \ \ Xiaolong Yang \thanks{School of Mathematics and Statistics, Central China Normal University, Wuhan 430079,  P. R. China.  Email: \texttt{yangxiaolong@mails.ccnu.edu.cn}.}}

\date{}


\smallbreak \maketitle

\begin{abstract}
We study the existence of normalized solutions to the following logarithmic Schr\"{o}dinger equation
\begin{equation*}\label{eqs01}
-\Delta u+\lambda u=\alpha u\log u^2+\mu|u|^{p-2}u, \ \ x\in\R^N,
\end{equation*}
under the mass constraint
\[
\int_{\R^N}u^2\mathrm{d}x=c^2,
\]
where $\alpha,\mu\in \R$, $N\ge 2$,  $p>2$, $c>0$ is a constant, and $\lambda\!\in\!\R$ appears as Lagrange multiplier.
Under different assumptions on $\alpha,\mu,p$ and $c$, we prove the existence of ground state solution and excited state solution.
The asymptotic behavior of the ground state solution as $\mu\to 0$
is also investigated. Our results including the case $\alpha<0$ or $\mu<0$, which is less studied in the literature.
%
\end{abstract}

\vskip 0.12truein


\vskip 0.08truein

\noindent{\text{2010 AMS Subject Classification:}}
35J91 (35J20 35B33 35B40).

\vskip 0.2truein

\section{Introduction}
In this paper, for prescribed $c>0$, we look for solutions $(\lambda,u)\in \R\times H^1(\R^N)$ satisfying
\begin{equation}\label{eq1.1}
-\Delta u+\lambda u=\alpha u\log u^2+\mu|u|^{p-2}u, \ \ x\in\R^N,
\end{equation}
and
\begin{equation}\label{eq1.11}
\int_{\R^N}u^2\mathrm{d}x=c^2,
\end{equation}
where $\alpha,\mu\in\R$ and $2<p\le 2^*=\frac{2N}{N-2}$ ($2^*=\infty$ if $N=1,2$).
Equation \eqref{eq1.1} is closely related to the following
time-dependent Schr\"{o}dinger equation
\begin{equation}\label{eq1.2}
-i\partial_t\Phi-\Delta \Phi-\alpha \Phi \log|\Phi|^2-\mu|\Phi|^{p-2}\Phi=0\quad\ \text{in}\ \ \R^N.\\
\end{equation}
Equation \eqref{eq1.2} admit plenty of applications related to quantum mechanics,
quantum optics, nuclear physics, transport and diffusion phenomena, open quantum
systems, effective quantum gravity, theory of superfluidity and Bose-Einstein condensation,
see \cite{AC,CP,CT1,SS,T,WZ} and the references therein.

\vskip1mm

Normalized solutions for equation \eqref{eq1.1} seem to be particularly meaningful from the physical
viewpoint, as pointed out in \cite{Soave1,SN}, these solutions often offer a good insight
of the dynamical properties of the stationary solutions for the nonlinear Schr\"odinger
equation \eqref{eq1.2}, such as stability and instability \cite{BC,CL}.

\vskip1mm

In the past three decades, the existence and multiplicity of normalized solutions for Schr\"odinger equation has been investigated by many authors, see \cite{BS1,BS,CT,CJ,BL1,BL,JJ,JJL,JTT,JL,JLS1,JLS2,MS21,Soave1,SN,WW,GS1,GLW,BMR}. In particular, if we replace the term $\alpha u\log u^2$ by $|u|^{q-2}u$,
then equation \eqref{eq1.1} turns into
\begin{equation}\label{eqs-p-q}
-\Delta u+\lambda u= |u|^{q-2}u+\mu|u|^{p-2}u, \ \  \ x\in \R^N,
\end{equation}
where $\mu\in\R$. If $2<p\le 2+\frac{4}{N}\le q\le 2^*$, Soave et al. \cite{Soave1,SN,WW,JJL,JTT} studied the existence and multiplicity of normalized solutions for equation \eqref{eqs-p-q}. Thereafter, normalized solutions to elliptic PDEs have attracted much
attention of researchers e.g. \cite{JL,JLS1,JLS2,MS21,BMR}.
\vskip1mm

In the case $\mu=0$, equation \eqref{eq1.1} is the so-called logarithmic Schr\"{o}dinger equation.
Recently,
there has been increasing interest in studying logarithmic Schr\"{o}dinger equation,
specially on the existence
of positive solutions, multiple solutions, ground states and semiclassical states, see for
examples, \cite{AJ,CG,ITWZ,SW,ZW,T,TZ}. In particular, Cazenave \cite{CT}
investigated the existence and multiplicity of normalized solutions
for equation \eqref{eq1.1}--\eqref{eq1.11} with $\alpha=1$ and $\mu=0$.


%
%

\vskip1mm

Inspired by the above results, an interesting question is that whether one can find normalized solutions for equation \eqref{eq1.1}--\eqref{eq1.11}. The goal of the present paper
is to give an affirmative answer. Although equation \eqref{eq1.1} has, at least formally, a variational structure related to the energy functional
\begin{equation*}
\begin{aligned}
I_{\lambda}(u)&=\frac{1}{2}\int_{\R^N}|\nabla u|^2+(\lambda+\alpha) u^2\mathrm{d}x-\frac{\alpha}{2}\int_{\R^N}u^2\log u^2\mathrm{d}x-\frac{\mu}{p}\int_{\R^N} |u|^p\mathrm{d}x,\\
\end{aligned}
\end{equation*}
this energy functional is not well defined on the natural Sobolev space $H^1(\R^N)$. In \cite{CT}, the idea is to work on the Banach space
\begin{equation}\label{space W}
W=\Big\{u\in H^1(\R^N)~|~\int_{\R^N}u^2|\log u^2|\mathrm{d}x<\infty \Big\},
\end{equation}
which is equipped the norm given by
\begin{equation*}
\|u\|_{W}=\|u\|_{H^1(\R^N)}+\inf\Big\{k>0~|~\int_{\R^N}A(k^{-1}|u|)\mathrm{d}x\le 1 \Big\},
\end{equation*}
where $A(s)=-s^2\log s^2$ on $[0,e^{-3}]$ and $A(s)=3s^2+4e^{-3}s-e^{-6}$ on $[e^{-3},\infty)$. In fact, by \cite[Proposition 2.7]{CT},
$I_{\lambda}:W\to \R$ is well defined and $\mathcal{C}^1$ smooth.


\vskip0.8mm

Solutions to equations \eqref{eq1.1}-\eqref{eq1.11} can be obtained by finding critical points of the energy functional
\begin{equation*}
I(u)=\frac{1}{2}\int_{\R^N}|\nabla u|^2+\alpha u^2\mathrm{d}x-\frac{\alpha}{2}\int_{\R^N} u^2\log u^2\mathrm{d}x-\frac{\mu}{p}\int_{\R^N}|u|^p\mathrm{d}x
\end{equation*}
under the constraint
\begin{equation*}
S(c):=\Big\{u\in W ~\big|~\|u\|_2=c \Big\}.
\end{equation*}

\vskip1.5mm

Before introducing our main results, we give some definitions (see also \cite{BeJe}).
\begin{defn}
We say that $u_0$ is a ground state solution for equation \eqref{eq1.1}--\eqref{eq1.11} on $S(c)$ if
\begin{equation*}
dI \vert_{S(c)}(u_0)=0 \quad \hbox { and } \quad I(u_0)=\inf \Big\{I(u) \hbox{ s.t. } dI\vert_{S(c)}(u)=0\ \hbox{and}\ u\in S(c)\Big\}.
\end{equation*}
We say that $v_0$ is an excited state solution for \eqref{eq1.1}--\eqref{eq1.11} on $S(c)$ if
\begin{equation*}
dI\vert_{S(c)}(v_0)=0 \quad \text { and } \quad I(v_0)>\inf \Big\{I(u) \hbox{ s.t. } dI\vert_{S(c)}(u)=0\ \text{and}\ u\in S(c)\Big\}.
\end{equation*}
\end{defn}

\vskip1.5mm

For $2<p<2^*$, the Gagliardo-Nirenberg inequality (see \cite{WeM}) is
\begin{equation}\label{b2}
\|u\|_p\leq  C(N,p)\|\nabla u\|_2^{\gamma_p}\|u\|_2^{1-\gamma_p},\quad \text{for each}\ \ u\in H^1(\R^N),
\end{equation}
where $C(N,p)$ is the best constant in the Gagliardo-Nirenberg inequality, $\gamma_p=\frac{N(p-2)}{2p}$.

\vskip1.5mm

We shall mainly concentrate here on the cases $\alpha=\pm 1$ and $\mu\in \R$. For $\alpha=0$, we refer the reader to \cite{Soave1,SN}.
Our main results of the paper can be stated as follows. We first consider the existence of global minima on $S(c)$.
Define
\begin{equation*}
m(c):=\inf_{u\in S(c)}I(u).
\end{equation*}

\begin{thm}\label{th1.1}
Let $\alpha=1$. Suppose $N\ge 2$, $c>0$ and that one of the following three conditions holds
\begin{equation*}
\begin{aligned}
&(i)~ \mu\le 0 \ \text{and} \ p>2; \qquad (ii)~\mu>0\ \text{and} \ 2<p<2+\frac{4}{N};\\
&(iii)~\mu>0,~p=\bar{p}:=2+\frac{4}{N},\ \ \text{and}\ \ c<\Big(\frac{N+2}{\mu N}\Big)^{\frac{N}{4}}\Big(\frac{1}{C(N,\bar{p})}\Big)^{\frac{N+2}{2}};\\
\end{aligned}
\end{equation*}
where $C(N,\bar{p})$ is the best constant of \eqref{b2} with $p=\bar{p}$, then
the infimum $m(c)$ admits a minimizer $\tilde{u}\in S(c)$,
which is positive, radially symmetric and decreasing in $r=|x|$.
Moreover, $\tilde{u}$ is a ground state of \eqref{eq1.1}-\eqref{eq1.11}.
\end{thm}

If $2+\frac{4}{N}< p \le 2^*$, we can easily verify $\inf\limits_{u\in S(c)}I(u)=-\infty$.  In spirit of \cite{BS}, we introduce the following Pohozaev manifold
\begin{equation*}\label{c1}
\mathcal{P}_{c}=\big\{u \in S(c) ~|~ P(u)=0\big\},
\end{equation*}
where
\begin{equation}\label{b3}
\begin{aligned}
P(u) =:\|\nabla u\|^2_2-\mu\gamma_p\|u\|^p_p-\frac{N}{2}\alpha c^2,
\end{aligned}
\end{equation}
with $\gamma_p=\frac{N(p-2)}{2p}$.
If $u \in W$ is a weak solution of \eqref{eq1.1}, then the Pohozaev identity $P(u)=0$ holds.
For $u \in S(c)$ and $s\in \mathbb{R}$, we define
\begin{equation*}
 s\star u(x) :=e^{\frac{N}{2}s} u\left(e^s x\right),~~~~\mbox{for}~\mbox{a.e.}~x \in \mathbb{R}^{N}.
\end{equation*}
The Pohozaev manifold $\mathcal{P}_{c}$ is quite closely related to the fiber map
\begin{equation*}
\begin{aligned}
\Psi_{u}(s):=I(s\star u)&=\frac{e^{2s}}{2}\int_{\R^N}|\nabla u|^2\mathrm{d}x +\frac{\alpha}{2}\int_{\R^N}u^2\mathrm{d}x-\alpha\frac{Ns}{2}\int_{\R^N} u^2\mathrm{d}x \\
&\quad -\frac{\alpha}{2}\int_{\R^N}u^2\log u^2\mathrm{d}x -\frac{e^{p\gamma_ps}}{p}\mu\int_{\R^N}|u|^p\mathrm{d}x.\\
\end{aligned}
\end{equation*}
For $u \in S(c)$ and $s \in \R$, we have
\begin{equation*} \label{c2}
\Psi'_{u}(s)=P\big(s \star u\big),
\end{equation*}
where $P$ is defined by \eqref{b3}.
We shall see that the critical points of $\Psi_{u}(s)$ allow us to project a function on $\mathcal{ P}_{c}$. Thus, the monotonicity and convexity properties of $\Psi_{u}(s)$ strongly affect the structure of $\mathcal{ P}_{c}$.
Generally, $\mathcal{P}_{c}$ can be divided into the disjoint union $\mathcal{ P}_{c}=\mathcal{ P}_{c}^+\cup \mathcal{ P}_{c}^0\cup \mathcal{ P}_{c}^-$, where
\begin{equation*}\label{c41}
\begin{aligned}
	\mathcal{ P}_{c}^+&:=\big\{u\in S(c)~|~\Psi_{u}'(0)=0,\  \Psi_{u}''(0)>0\big\},\\
	\mathcal{ P}_{c}^0&:=\big\{u\in S(c)~|~\Psi_{u}'(0)=0,\  \Psi_{u}''(0)=0\big\},\\
	\mathcal{ P}_{c}^-&:=\big\{u\in S(c)~|~\Psi_{u}'(0)=0,\  \Psi_{u}''(0)<0\big\}.
\end{aligned}
\end{equation*}
Define
\begin{equation*}
m^{+}(c):=\inf_{u\in \mathcal{P}^{+}_{c}} I(u),\qquad \ m^{-}(c):=\inf_{u\in \mathcal{P}^{-}_{c}} I(u)
\end{equation*}
and
$$c_0:=\Big[\frac{p2^{\frac{p\gamma_p}{2}}}{\mu \big(p\gamma_p\big)^{\frac{p\gamma_p+2}{2}}}\Big(\frac{p\gamma_p-2}{N}\Big)^{\frac{p\gamma_p-2}{2}}C^{-p}(N,p)\Big]^{\frac{1}{p-2}}.$$

If $\alpha=1$, $\mu>0$, $2+\frac{4}{N}<p<2^*$, then the functional
$I|_{S(c)}$ present a convex-concave geometry structure for $c>0$ small.
We thus show that \eqref{eq1.1}-\eqref{eq1.11} admits two solutions $u_c^+$ and $u_c^{-}$, which can be characterized respectively as local minima and mountain pass critical point of $I$ restricted to $S(c)$.

\begin{thm}\label{th1.2}
Let $\alpha=1$, $\mu>0$, $2+\frac{4}{N}<p<2^*$, $N\ge 2$. If $0<c<c_0$, then

(1)\,\, $m^{+}(c)$ is achieved by some $u_c^{+}$;
(2)\,\, $m^{-}(c)$ is achieved by some $u_c^{-}$.\\
Moreover, $u_c^{+}$, $u_c^{-}$ are positive, radially symmetric and decreasing in $r=|x|$. In addition there exist $\lambda_c^+$, $\lambda_c^-\in\R$ such that $(\lambda_c^+,u^+_c)$ and $(\lambda_c^-,u^-_c)$ are solutions to equation \eqref{eq1.1}-\eqref{eq1.11}, and $u^+_c$ is a ground state solution, $u^-_c$ is a excited  solution.
\end{thm}

Next, we study the existence of normalized solutions to equation \eqref{eq1.1}--\eqref{eq1.11} with Sobolev critical exponent.

\begin{thm}\label{th1.3}
Let $\alpha=1$, $\mu>0$, $p=2^*$, $N\ge 3$. If $0<c<c_0$, then

(1)\,\, $m^{+}(c)$ is achieved by some $u_c^{+}$;
(2)\,\, $m^{-}(c)$ is achieved by some $u_c^{-}$.\\
Moreover, $u_c^{+}$, $u_c^{-}$ are positive, radially symmetric and decreasing in $r=|x|$. In addition there exist $\lambda_c^+$, $\lambda_c^-\in\R$ such that $(\lambda_c^+,u^+_c)$ and $(\lambda_c^-,u^-_c)$ are solutions to equation \eqref{eq1.1}-\eqref{eq1.11}, and $u^+_c$ is a ground state solution, $u^-_c$ is a excited  solution.
\end{thm}

\begin{Rem}

(1) For $2+\frac{4}{N}<p\le 2^*$ and $0<c\le c_0$, we first study the following local minimization problem
\begin{equation*}
m_c:=\inf_{u\in V_{k_0}} I(u),\quad \text{where} \quad V_{k_0}=\big\{u\in S(c)~|~\|\nabla u\|^2_2<k_0 \big\}, \ k_0=\frac{p\gamma_pNc^2}{2(p\gamma_p-2)}.
\end{equation*}
In order to find more than one solution, we decompose the natural constraint $\mathcal{P}_c$ into three disjoint subsets $\mathcal{P}^+_c$,
$\mathcal{P}^0_c$ and $\mathcal{P}^-_c$. If $0<c<c_0$, then $u_c^+$ is a global minimizer of $I$ restricted on $\mathcal{P}_c$, which can characterized as a local minimizer of $I$ on the set $V_{k_0}$. Moreover, the second solution $u_c^-$ corresponds to a critical point of mountain-pass type for $I$ on $S(c)$. It is worth mentioning that the existence of such two critical points on $S(c)$ for Schr\"{o}dinger equation has been studied in \cite{Soave1,CJ,SN,JL,WW}.

\vskip1mm

(2) If $p=2^*$, we have $\gamma_p=1$ and $C(N,p)=\mathcal{S}^{-\frac{1}{2}}$, then $c_0$ can be written as
$c_0:=\Big[\frac{N^2-2N}{4\mu}\Big(\frac{4\mathcal{S}}{N^2}\Big)^{\frac{2^*}{2}}\Big]^{\frac{1}{2^*-2}}$, where $\mathcal{S}$ is the optimal constant in the Sobolev embedding. Compared to the works \cite{CJ,JL}, we can even prove the existence of at least one positive solution for $c=c_0$, see Lemma \ref{lem5.2}.

\vskip1mm
(3) The condition $0<c<c_0$ in Theorems \ref{th1.2} and \ref{th1.3} not only ensures that the corresponding energy functional $I$ admits
a convex-concave geometry, but also guarantees that the Pohozaev manifold $\mathcal{P}_{c}$ is a natural constraint, on which
the critical points of $I$ are indeed normalized solutions for \eqref{eq1.1}--\eqref{eq1.11}. For $p=2^*$, we drive a better energy estimate on the associated mountain pass level 
\[
m^-(c)<\frac{\mu^{-\frac{N-2}{2}}}{N}\mathcal{S}^{\frac{N}{2}}+m^+(c),
\]
i.e., the mountain pass energy level is less than the usual critical threshold plus the ground state energy,
and thus ensures the compactness of Palais-Smale sequence and a mountain pass type solution follows.
\end{Rem}

\vskip1mm

Now, we study the asymptotic behavior of ground states as $\mu\to 0$.
\begin{thm}\label{th1.4}
Assume $(\lambda_{\mu},u^+_{\mu})$ is the normalized ground state solution of \eqref{eq1.1}-\eqref{eq1.11} obtained by Theorem \ref{th1.1} (Theorem \ref{th1.2} or Theorem \ref{th1.3}), then, up to a subsequence, we have
\begin{equation*}
u^+_{\mu}\to w_0\quad \text{strongly~in}\ \ W,
\end{equation*}
and $\lambda_{\mu}\to \lambda_0$ as $\mu \to 0$, where $(\lambda_0,w_0)$ is a normalized ground state solution for \eqref{eq1.1}-\eqref{eq1.11} with $\mu=0$.
\end{thm}

Next, we consider the case of $\alpha=-1$. The term $-u\log|u|^2$ makes the
geometry structure of $I|_{S(c)}$ much more complex,
we follow the ideas introduced by Cingolani \& Jeanjean \cite{CJ} and Jeanjean \& Le \cite{JL} to catch the normalized solutions.
The authors in \cite{CJ,JL} studied the existence of normalized solutions for the following Schr\"odinger-Poisson system
\begin{equation*}
\begin{cases}
-\Delta u +\lambda u+\gamma \phi u=a|u|^{p-2}u\quad \text{in}\ \R^d,\\
-\Delta \phi=u^2 \quad \text{in}\ \R^d,
\end{cases}
\end{equation*}
for $d=2$ or $d=3$ respectively.
Under different assumptions on $\gamma$, $a\in \R$ and $p$, they proved several existence and multiplicity results.

\begin{thm}\label{th1.51}
\noindent \textup{(i)} If $p:=\bar{p}=2+\frac{4}{N}$, $\mu>0$ and $c<\Big(\frac{N+2}{\mu N}\Big)^{\frac{N}{4}}\big(C(N,\bar{p})\big)^{-\frac{N+2}{2}}$,
then $I$ does not have critical points on $S(c)$.

\noindent \textup{(ii)} If $\mu\le 0$, for any $p>2$ and $c>0$, then $I$ does not have critical points on $S(c)$.

\noindent \textup{(iii)} Let $p>2$, $N\ge 2$ and $c>0$, then equation \eqref{eq1.1}-\eqref{eq1.11}
has no positive radial solution in $W\cap H^1_r(\R^N)$ for all $\mu \in \R$.
\end{thm}

Define
\begin{equation*}
D:=\Big(\frac{Np\gamma_p}{2(2-p\gamma_p)}\Big)^{\frac{2-p\gamma_p}{2(p-2)}}\Big(\frac{\mu p\gamma^2_p}{2}C^p(N,p)\Big)^{-\frac{1}{p-2}}.
\end{equation*}
By Lemma \ref{lem5.2}, we obtain
\begin{equation*}
\inf_{u\in \mathcal{P}_c}I(u)=-\infty\quad \text{if}\ c\ge D.
\end{equation*}
However
\begin{equation*}
\sup_{u\in \mathcal{P}_c} I(u)<+\infty\quad \text{if}\  c \ge D.
\end{equation*}
We are able to prove that the existence of a global maximal on $\mathcal{P}_c$.
\begin{thm}\label{th1.5}
Let $\alpha=-1$. If $p\in(2,2+\frac{8}{N(N+2)})\bigcup (2+\frac{8}{N(N+2)},2+\frac{4}{N})$, $N\ge 2$, $\mu>0$ and $c=D$, there exists $\bar{u}\in S(c)$ such that
\begin{equation*}
I(\bar{u})=\sup_{u\in \mathcal{P}_c\cap H^1_r(\R^N)} I(u).
\end{equation*}
Moreover, $\bar{u}$ is a critical point of $I$ restricted to $S(c)$, and $\bar{u}$ is non-positive and radially symmetric.
\end{thm}

\begin{Rem}
(1) If $\alpha<0$, for fixed frequency case, the study of \eqref{eq1.1} is still an open field of investigation.
Fortunately, with the help of prescribed $L^2$-norm and Pohozaev manifold, we get partial existence and non-existence results.

\vskip1mm

(2) By Lemma \ref{lem7.1}, $\mathcal{P}_c$ is not empty if and only if $c\ge D$. In addition, if $0<c\le D$, $\mathcal{P}_c$ is a smooth manifold of codimension 2 in $W\cap H^1_r(\R^N)$ (see Lemma \ref{lem7.5}). Therefore, for $c=D$, we prove the existence of one global maximizer of $I$ on $\mathcal{P}_c$. In particular, in the mass-subcritical case $2<p<2+\frac{4}{N}$,
it could be natural to expect that there exists a second radial critical point on $S(c)$. However, we can not
ensure that the corresponding energy functional $I$ admits a concave-convex geometry structure on $S(c)$, and can not guarantees that the Pohozaev manifold $\mathcal{P}_{c}$ is a natural constraint.
\end{Rem}

This paper is organized as follows. In Section 2, we give some preliminary results. In Section 3, we study the existence of the associated global minimizer and prove Theorem \ref{th1.1}.
In Section 4, we mainly focus on mass-supercritical and Sobolev subcritical case.
In Section 5, mass-supercritical and Sobolev critical case is investigated, two pairs normalized solutions are obtained.
While in Section 6, we give a precise asymptotic behavior of the normalized ground state solutions as $\mu\to 0$.
Finally, in Section 7, we prove some existence and non-existence results of \eqref{eq1.1}-\eqref{eq1.11} in the case of $\alpha<0$.


\vskip2mm

{\bf Notations.} In the paper, we use the following notations. $L^p=L^p(\R^N)$ with norm $\|u\|_{L^p(\R^N)}=\|u\|_p$,
$H^1(\R^N)$ is the usual Sobolev space with norm $\|u\|_{H^1(\R^N)}=\big(\int_{\R^N}|\nabla u|^2+u^2\mathrm{d}x\big)^{\frac{1}{2}}$. $C, C_i, i=1,2,\cdots$ and so on denote universal positive constants, which we need not to specify, and which may vary from line to line.

\section{Preliminaries}
In this section, we give some preliminary results.

Similar as \cite{CT}, we define
\begin{equation*}
A(s):=
\begin{cases}
-s^2\log s^2,\qquad & \text{if} ~~0\le s\le e^{-3},\\
3s^2+4e^{-3}s-e^{-6},\quad &\text{if} ~~ e^{-3}\le s,\\
\end{cases}
\end{equation*}
and
\begin{equation*}
B(s):=s^2\log s^2+A(s).
\end{equation*}
From Lemma 1.2 in \cite{CT}, $A$ is a positive convex increasing function. Moreover, for every $s>0$, there exists $K_q>0$ such that
\begin{equation}\label{inequality B}
B(s)\le K_q s^q \quad \text{for~all}~q\in \big(2,2+\frac{4}{N}\big).
\end{equation}
Denote
\begin{equation*}
V:=\Big\{ u\in L^1_{loc}(\R^N)~|~A(|u|)\in L^1(\R^N)  \Big\}.
\end{equation*}

\begin{lem}\label{lem2.1}\cite[Lemma 2.1]{CT}
$(i)$ $V$ is the Orlicz space associated to $A$. $V$ equipped with the norm $\|\cdot\|_V$ defined by
\begin{equation*}
\|u\|_V:=\inf\Big\{k>0~|~\int_{\R^N}A(k^{-1}|u|)\mathrm{d}x\le 1\Big\}
\end{equation*}
is a reflexive Banach space.
\vskip1mm
$(ii)$ For any $u\in V$, $\inf\big\{\|u\|_V,\|u\|^2_{V}\big\}\le \int_{\R^N}A(|u|)\mathrm{d}x\le \sup\big\{\|u\|_{V},\|u\|^2_{V} \big\}.$

\vskip1mm
$(iii)$ If $u_n\to u$ a.e. in $\R^N$ and $\int_{\R^N}A(|u_n|)\mathrm{d}x\to \int_{\R^N}A(|u|)\mathrm{d}x<\infty$, then $\|u_n-u\|_V\to 0$ as $n\to \infty$.
\end{lem}


\vskip1.2mm

Define
\begin{equation*}
W_r:=W\cap H^1_r(\R^N),
\end{equation*}
where $W$ is defined by \eqref{space W}. Obviously, $W$ is a reflexive Banach space, since $H^1(\R^N)$ and $V$ are both reflexive Banach spaces.

\vskip1.2mm

\begin{lem}\label{lem2.2}\cite[Proposition 2.7, Proposition 3.1]{CT} The following facts holds

$(i)$ $I\in \mathcal{C}^1(W,\R)$ and for any $u\in W$ one has $DI(u)=Lu$, where
$$Lu=-\Delta u-\alpha u\log u^2-\mu |u|^{p-2}u.$$
\vskip1mm
$(ii)$ The embedding from $W_r$ to $L^2(\R^N)$ is compact.
\end{lem}

\begin{lem}\cite[Lemma 3.6]{BS1}\label{lem2.3}
For $u\in S(c)$ and $t\in\R$,
the map $T_{u}S(c)\to T_{t\star u}S(c)$ defined by $\psi\to t\star\psi$ is a linear isomorphism, where
\[
T_{u}S(c)=\Big\{v\in H^1(\R^N)~\big|~ \int_{\R^N}uv\mathrm{d}x=0\Big\}.
\]
\end{lem}





\section{The proof of Theorem \ref{th1.1}}
In this section, we prove Theorem \ref{th1.1}.

Define
\begin{equation*}\label{a1}
m(c):=\inf_{u\in S(c)}I(u)\ \ \ \text{and}\ \ \ m_r(c):=\inf\limits_{u\in S(c)\cap H^1_r(\R^N)}I(u).
\end{equation*}

\begin{lem}\label{lem3.1}
Under the assumption of Theorem \ref{th1.1}, then $m(c)>-\infty$
and $m(c)=m_r(c)$.
\end{lem}

\noindent{\bf Proof.}
We first prove $m(c)>-\infty$. For case $(i)$: $\mu\le 0 \ \text{and} \ p>2$. Let $2<q<2+\frac{4}{N}$, for each $u\in W$, by using \eqref{inequality B}, we can deduce that
\begin{equation}\label{a2}
\begin{aligned}
I(u)&\ge \frac{1}{2}\int_{\R^N}|\nabla u|^2+u^2\mathrm{d}x+\frac{1}{2}\int_{\R^N}A(|u|)\mathrm{d}x-\frac{1}{2}\int_{\R^N}B(|u|)\mathrm{d}x\\
&\ge\frac{1}{2}\int_{\R^N}|\nabla u|^2+u^2\mathrm{d}x-\frac{1}{2}K_q\int_{\R^N}|u|^q\mathrm{d}x\\
&\ge\frac{1}{2}\|\nabla u\|^2_2-\frac{1}{2}C^q(N,q)K_qc^{q(1-\gamma_q)}\|\nabla u\|^{q\gamma_q}_2.
\end{aligned}
\end{equation}
For case $(ii)$: $\mu>0$ and $2<p<2+\frac{4}{N}$. Similar as \eqref{a2}, for each $u\in W$, we have
\begin{equation}\label{a3}
\begin{aligned}
I(u)&\ge\frac{1}{2}\int_{\R^N}|\nabla u|^2+u^2\mathrm{d}x
+\frac{1}{2}\int_{\R^N}A(|u|)-B(|u|)\mathrm{d}x-\frac{\mu}{p}\int_{\R^N}|u|^p\mathrm{d}x\\
&\ge\frac{1}{2}\|\nabla u\|^2_2-\frac{1}{2}C^q(N,q)K_qc^{q(1-\gamma_q)}\|\nabla u\|^{q\gamma_q}_2-\frac{\mu}{p}C^p(N,p)c^{p(1-\gamma_p)}\|\nabla u\|^{p\gamma_p}_2.
\end{aligned}
\end{equation}
For case $(iii)$: $\mu>0, p=\bar{p}$, and $c<\left(\mu^{-1}\frac{N+2}{N}\right)^{\frac{N}{4}}\big(C(N,\bar{p})\big)^{-\frac{N+2}{2}}$. It follows from \eqref{a3} that
\begin{equation}\label{a3-1}
\begin{aligned}
I(u)
&\ge\Big(\frac{1}{2}-C^{\bar{p}}(N,\bar{p})\frac{\mu c^{\frac{4}{N}}}{\bar{p}}\Big)\|\nabla u\|^2_2-C^q(N,q)K_qc^{q(1-\gamma_q)}\|\nabla u\|^{q\gamma_q}_2.
\end{aligned}
\end{equation}
Since $q\gamma_q$, $p\gamma_p<2$ for $p,q\in [2,2+\frac{4}{N})$,
we then conclude that, under the assumption of Theorem \ref{th1.1}, $I$ is coercive on $S(c)$. Therefore, $m(c)>-\infty$.

\vskip 0.08truein

On the other hand, one can easily conclude $m(c)\leq m_r(c)$, so we only need to prove $m(c)\geq m_r(c)$.

Assume that $\{u_n\}\subset S(c)$ is a minimizing sequence for $m(c)$. Denote $u_n^*$ be the symmetric decreasing rearrangement of
$u_n$. Thus, by \cite[(iv)--(v) of Chapter 3.3, and Lemma 7.17]{LL}, one has
\[
\int_{\R^N}|\nabla u_n^*|^2 dx\leq \int_{\R^N}|\nabla u_n|^2 dx,\ \
\int_{\R^N}|u_n^*|^r dx= \int_{\R^N}|u_n|^r dx \ \ \text{for}\ r\in [2,2^*].
\]
Recall that $A$, $B$ are positive, convex, increasing functions on $(0,+\infty)$, by using \cite[(v) of Chapter 3.3]{LL}, we get
\[
\int_{\R^N}A(u_n^*) dx= \int_{\R^N}A(u_n) dx,\ \
\int_{\R^N}B(u_n^*) dx= \int_{\R^N}B(u_n) dx,
\]
which implies
\[
\int_{\R^N}|u_n^*|^2 \log |u_n^*|^2 dx= \int_{\R^N}|u_n|^2 \log |u_n|^2 dx.
\]
Therefore
\[
m_r(c)=\inf\limits_{u\in S(c)\cap H^1_r(\R^N)}I(u)\leq \inf_{u\in S(c)}I(u) =m(c).
\]
\qed

\begin{lem}\label{lem3.2}
Under the assumption of Theorem \ref{th1.1}, then the infimum $m(c)$ is achieved by some $u\in S(c)$, which is a positive,
radially symmetric and decreasing in $r=|x|$.
\end{lem}

\noindent{\bf Proof.}
Assume that $\{u_n\}\subset S(c)\cap H^1_r(\R^N)$ be a minimizing sequence for $m(c)$. From \eqref{a2}--\eqref{a3-1},
we deduce that $\{u_n\}$ is bounded in $H^1(\R^N)$, and $\{\int_{\R^N}A(|u_n|)\mathrm{d}x\}$ is bounded. By Lemma \ref{lem2.1}, we conclude $\{u_n\}$ is bounded in $W_r$. Therefore, $u_n\rightharpoonup u$ weakly in $W_r$, it follow from Lemma \ref{lem2.2} that $u_n\to u$ strongly in $L^2(\R^N)$ and $u_n \to u$ a.e. in $\R^N$. If $\mu>0$, we deduce that
\begin{equation*}
B(|u_n|)\to B(|u|) \ \text{strongly~in} \ L^1(\R^N) \quad \text{and}\quad u_n\to u \ \text{strongly~in} \ L^p(\R^N).
\end{equation*}
Therefore
\begin{equation}\label{a4}
\begin{aligned}
m(c)&\le I(u)=\frac{1}{2}\int_{\R^N}|\nabla u|^2+u^2\mathrm{d}x-\frac{1}{2}\int_{\R^N}u^2\log u^2\mathrm{d}x-\frac{\mu}{p}\int_{\R^N}|u|^p\mathrm{d}x\\
&\le \liminf_{n \to \infty}\Big(\frac{1}{2}\int_{\R^N}|\nabla u_n|^2
+u^2_n\mathrm{d}x+\frac{1}{2}\int_{\R^N}A(|u_n|)\mathrm{d}x\Big)\\
&\quad-\frac{1}{2}\int_{\R^N}B(|u|)\mathrm{d}x-\frac{\mu}{p}\int_{\R^N}|u|^p\mathrm{d}x\\
&\le \liminf_{n \to \infty}I(u_n)=m(c).
\end{aligned}
\end{equation}
Hence, $I(u)=m(c)$, $u_n\to u$ strongly in $H^1(\R^N)$ and $A(|u_n|)\to A(|u|)$ in $L^1(\R^N)$. It thus follows from Lemma \ref{lem2.1}
that $u_n\to u$ strongly in $V$ as $n\to \infty$. Therefore, $u_n\to u$ strongly in $W_r$.

If $\mu\le 0$, we apply the same argument. By the weak lower semi-continuity, we have
\[
\liminf\limits_{n\to\infty}\Big(-\frac{\mu}{p}\int_{\R^N}|u|^p\mathrm{d}x\Big)
\le-\frac{\mu}{p}\int_{\R^N}|u_n|^p\mathrm{d}x.
 \]
Similar as \eqref{a4}, we obtain $u_n\to u$ strongly in $W_r$.
Since $u$ is a nonnegative nontrivial weak solution of \eqref{eq1.1}-\eqref{eq1.11}. Moreover, by elliptic regularity theory, we obtain $u\in \mathcal{C}^2(\R^N)$.
For $a>0$ small enough, we have
\begin{equation*}
\Delta u=\lambda u-u\log u^2-\mu u^{p-1}\le \beta(u)  \quad \text{in}\ \{x\in \R^N~|~0<u(x)<a\},
\end{equation*}
where $\beta(s)=\lambda s-s\log s^2$ for $\mu>0$ and $\beta(s)=(\lambda-\mu)s-s\log s^2$ for $\mu\le 0$.
Since $\beta$ is continuous, nondecreasing for $s$ small, $\beta(0)=0$ and $\beta(\sqrt{e^{\lambda}})=0$ for $\mu>0$, $\beta(\sqrt{e^{\lambda-\mu}})=0$ for $\mu\le 0$,
by \cite[Theorem 1]{VJ}, we have that $u>0$.
\qed

\vskip 0.12truein

{\bf Proof of Theorem \ref{th1.1}.}
The proof follows directly from Lemma \ref{lem3.2}.


\section{Mass-supercritical and Sobolev subcritical case}

In this section, we deal with the mass supercritical case $2+\frac{4}{N}<p<2^*=\frac{2N}{(N-2)^+}$, $N\ge 2$,
$\alpha>0$ and $\mu>0$, the functional $I$ is unbounded from below on $S(c)$, it will not be possible to find a global minimizer.

\subsection{Existence of a local minima on $S(c)$}

For $\alpha=1$, $c>0$ and $2+\frac{4}{N}<p<2^*$, $N\ge 2$, define
\begin{equation}\label{z1}
k_0=\frac{p\gamma_pNc^2}{2(p\gamma_p-2)},  \qquad  c_0=\left[\frac{p2^{\frac{p\gamma_p}{2}}}{\mu \big(p\gamma_p\big)^{\frac{p\gamma_p+2}{2}}}\Big(\frac{p\gamma_p-2}{N}\Big)^{\frac{p\gamma_p-2}{2}}C^{-p}(N,p)\right]^{\frac{1}{p-2}}.
\end{equation}

\begin{lem}\label{lem4.1}
Suppose that $\mu>0, c>0$, and $2+\frac{4}{N}<p<2^*$. If $P(u)\le 0$ and $\|\nabla u\|^2_2=k_0$, then $c\ge c_0$. Moreover, if $P(u)\le 0$ and $c<c_0$, then $\|\nabla u\|^2_2\neq k_0$.
\end{lem}
\noindent{\bf Proof.}
Since $P(u)\le 0$, we have
\begin{equation*}
\|\nabla u\|^2_2\le \mu \gamma_p \|u\|^p_p+\frac{N}{2}c^2.
\end{equation*}
By using the Gagliardo-Nirenberg inequality \eqref{b2}, it follows that
\begin{equation*}
\|\nabla u\|^2_2\le \mu \gamma_p C^p(N,p)c^{p(1-\gamma_p)}\|\nabla u\|^{p\gamma_p}_2+\frac{N}{2}c^2.
\end{equation*}
If $\|\nabla u\|^2_2=k_0$, then
\begin{equation*}
\frac{N}{p\gamma_p-2}\le \mu \gamma_p C^p(N,p)\Big(\frac{Np\gamma_p}{2(p\gamma_p-2)}\Big)^{\frac{p\gamma_p}{2}}c^{p-2},
\end{equation*}
which follows that $c\ge c_0$. Therefore, we deduce that if $P(u)\le 0$ and $c<c_0$, then $\|\nabla u\|^2_2\neq k_0$. We obtain the conclusion of the lemma.
\qed

\vskip 0.12truein

Now, we define
\begin{equation}\label{y1}
V_{k_0}:=\Big\{u\in S(c)~|~ \|\nabla u\|^2_2<k_0\Big\}.
\end{equation}
where $k_0$ is defined by \eqref{z1}. Next, for each $0<c\le c_0$, we study the following minimization problem
\begin{equation}\label{z2}
m_c := \inf_{u \in V_{k_0}} I(u).
\end{equation}

\begin{lem}\label{lem4.2}
Let $\mu>0$ and $2+\frac{4}{N}<p<2^*$. If $0<c\le c_0$, then $m_c$ is achieved by some $u$, which
is a positive, radially symmetric critical point of $I$ on $S(c)$.
\end{lem}
\noindent{\bf Proof.}
Let $\{u_n\}$ be a minimizing sequence for $m_c$, similar as the proof of Lemma \ref{lem3.2}, we can deduce $u_n\to u$ strongly in $W_r$.
We only need to prove that $u\in V_{k_0}$. In fact, if $\|\nabla u\|^2_2=k_0$ and $0<c<c_0$, it directly follows from Lemma \ref{lem4.1} that $P(u)>0$. Therefore, there exist $t_0<0$ such that $t_0\star u\in V_{k_0}$, and
$I(t_0\star u)<I(u)=m_c.$ This is a contradiction. On the other hand, if $\|\nabla u\|^2_2=k_0$ and $c=c_0$, the discussion is divided into three cases.
 \vskip1mm
Case 1: If $P(u)<0$, similar to the proof in Lemma \ref{lem4.1},
we have $c>c_0$, which is a contradiction.

\vskip1mm
Case 2: If $P(u)=0$, we have
\begin{equation}\label{x71}
\|\nabla u\|^2_2\le \mu \gamma_p C^p(N,p)c^{p(1-\gamma_p)}\|\nabla u\|^{p\gamma_p}_2+\frac{N}{2}c^2,
\end{equation}
The equality in \eqref{x71} holds only for the best constant in the Gagliardo-Nirenberg inequality is achieved, from \cite{WeM}, $u$ satisfies
\begin{equation}\label{x72}
\|\nabla u\|^2_2=\|u\|^2_2.
\end{equation}
Combining $P(u)=0$ and \eqref{x72}, which contradicts with $N\ge 2$. Therefore,
\begin{equation*}
\|\nabla u\|^2_2< \mu \gamma_p C^p(N,p)c^{p(1-\gamma_p)}\|\nabla u\|^{p\gamma_p}_2+\frac{N}{2}c^2,
\end{equation*}
we deduce that $c>c_0$ because $\|\nabla u\|^2_2=k_0$. This is also a contradiction.

\vskip1mm
Case 3: If $P(u)>0$, then there is $t_1<0$ such that $t_1\star u\in V_{k_0}$, and
$I(t_1\star u)<I(u)=m_c$. This is a contradiction.
\vskip1mm
In summarization, $u\in V_{k_0}$ and $I(u)=m_c$. Therefore, by combining $\|\nabla u\|^2_2<k_0$ with \cite[Proposition A.1]{MS21}, we conclude $u$ is a positive, radially symmetric critical point of $I$ on $S(c)$.
\qed

\subsection{Multiplicity of solutions}

In this subsection, we are interested in the multiplicity of solutions. For any $0<c<c_0$, we prove that \eqref{eq1.1}-\eqref{eq1.11} admits
two solutions, which can be characterized respectively as a local minimizer
or a mountain pass critical point of the
energy functional $I$ restricted to $S(c)$.
We first study the structure of the Pohozaev manifold $\mathcal{P}_c$.
Recalling the decomposition of $\mathcal{P}_c=\mathcal{P}^+_{c}\cap \mathcal{P}^{0}_{c}\cap \mathcal{P}^-_{c}$.
\begin{lem}\label{lem4.3}
Let $\mu>0$ and $2+\frac{4}{N}<p<2^*$. If $0<c<c_0$, then $\mathcal{P}^0_{c}=\emptyset$, and the set $\mathcal{P}_{c}$ is a $\mathcal{C}^1$-submanifold of codimension 1 in $S(c)$.
\end{lem}
\noindent{\bf Proof.}
Assume by contradiction that there exists $u \in \mathcal{P}^0_{c}$ such that $P(u)=0$ and
\begin{equation*}
\Psi''_{u}(0)=2\int_{\mathbb{R}^{N}}|\nabla u|^{2}\mathrm{d}x
-p\gamma^2_p\mu\int_{\R^N}|u|^p\mathrm{d}x=0.
\end{equation*}
Let
\begin{equation*}
\begin{aligned}
f(t):&=t\Psi'_{u}(0)-\Psi''_{u}(0)\\
&=(t-2)\int_{\mathbb{R}^{N}}|\nabla u|^{2}\mathrm{d}x
-(t-p\gamma_p)\gamma_p\mu\int_{\R^N}|u|^{p}\mathrm{d}x-\frac{N}{2}tc^2,\\
\end{aligned}
\end{equation*}
we see that $f(t)=0,~\forall t\in \R$. Therefore, it follows from $f(p\gamma_p)=0$ that
\begin{equation}\label{z6}
(p\gamma_p-2)\int_{\mathbb{R}^{N}}|\nabla u|^{2}\mathrm{d}x=p\gamma_p\frac{N}{2}c^2.
\end{equation}
From $\eqref{z6}$, we have $\|\nabla u\|^{2}_2=k_0$. Since $f(2)=0$, and $\|\nabla u\|^2_2=k_0$, 
we get
\begin{equation*}
\begin{aligned}
N \le \mu (p\gamma_p-2)\gamma_p C^p(N,p)\Big(\frac{Np\gamma_p}{2(p\gamma_p-2)}\Big)^{\frac{p\gamma_p}{2}}c^{p-2}
\end{aligned}
\end{equation*}
which is a contradiction with $c<c_0$. This proves that $\mathcal{P}^0_{c}=\emptyset$.

\vskip1mm
We now prove that $\mathcal{P}_{c}$ is a smooth manifold of codimension 1 in $S(c)$. We know that $\mathcal{P}_{c}$ is defined by
$P(u)=0$ and $G(u)=0$ where
\begin{equation*}
G(u)=\|u\|^2_2-c^2.
\end{equation*}
Since $P$ and $G$ are of $\mathcal{C}^1$-class, the proof is complete provided we show that the differential
$(dP(u),dG(u)): W \to \R^2 $ is surjective,
for every \[ u\in G^{-1}(0)\cap P^{-1}(0). \]
If this is not true, $dP(u)$ has to be linearly dependent from $dG(u)$, i.e. there exists $\nu \in \mathbb{R}$  such that
\begin{equation*}
\int_{\R^N} 2\nabla u\nabla \varphi_1 +2\nu u\varphi_1\mathrm{d}x
=\int_{\R^N} p\gamma_p\mu |u|^{p-2}u\varphi_1\mathrm{d}x, \ \ \text{for each}\ \varphi_1 \in W.
\end{equation*}
Therefore, $u$ satisfies
\begin{equation*}
-2\Delta u +2\nu u=p\gamma_p\mu|u|^{p-2}u.
\end{equation*}
The Pohozaev identity for the above equation is
\begin{align*}
2\|\nabla u\|^2_2=p\gamma^2_p\mu\|u\|^p_p,
\end{align*}
which contradicts to $\mathcal{P}^0_{c}=\emptyset$.
\qed

\begin{lem}\label{lem4.4}
If $2+\frac{4}{N}<p< 2^*$, $N\ge 2$ and $0<c<c_0$, for each $u\in S(c)$, then the function
$h(s):=\Psi_{u}(s)$ has exactly two critical points $s_u$, $t_u$ with $s_{u}\!<\!t_{u}$. Moreover,

\begin{itemize}
\item  [$(i)$] $s_{u}\star u \!\in\!\mathcal{P}^{+}_{c}$ and $s_u$ is a strict local minimum point for $h(s)$; $t_{u}\star u \!\in\!\mathcal{P}^{-}_{c}$ and $t_u$ is a strict local maximum point for $h(s)$;

\item  [$(ii)$] The maps $u \mapsto s_{u} \in \mathbb{R}$
and $u \mapsto t_{u} \in \mathbb{R}$ are of class $\mathcal{C}^1$.
\end{itemize}
\end{lem}

\noindent{\bf Proof.}
$(i)$ For each fixed $u\in S(c)$ and $\mu>0$, there exists $\bar{t}=\frac{1}{p\gamma_p-2}\log\Big(\frac{2\|\nabla u\|^2_2}{p\gamma^2_p\mu\|u\|^p_p}\Big)$ such that
\begin{equation*}
2e^{2\bar{t}}\|\nabla u\|^2_2=p\gamma^2_pe^{p\gamma_p\bar{t}}\mu\|u\|^p_p,
\end{equation*}
i.e.
$$2\|\bar{t}\star(\nabla u)\|^2_2=p\gamma^2_p\mu\|\bar{t}\star u\|^p_p.$$
It follows that
\begin{equation}\label{z7}
2\|t\star(\nabla u)\|^2_2>p\gamma^2_p\mu\|t\star u\|^p_p\quad \text{for~all}~-\infty<t<\bar{t}
\end{equation}
and
\begin{equation*}
2\|t\star(\nabla u)\|^2_2<p\gamma^2_p\mu\|t\star u\|^p_p\quad \text{for~all}~t>\bar{t}.
\end{equation*}
From \eqref{z7}, for $-\infty<t<\bar{t}$, we have
\begin{equation}\label{z8}
\begin{aligned}
\Psi'_{u}(t)&=\|t\star (\nabla u)\|^2_2-\mu \gamma_p\|t\star u\|^p_p-\frac{N}{2}c^2\\
&>\|t\star (\nabla u)\|^2_2-\frac{2}{p\gamma_p}\|t\star(\nabla u)\|^2_2-\frac{N}{2}c^2\\
&=\frac{p\gamma_p-2}{p\gamma_p}\|t\star (\nabla u)\|^2_2-\frac{N}{2}c^2.
\end{aligned}
\end{equation}
Now, we claim that if $c<c_0$, then
\begin{equation*}
\|\bar{t}\star (\nabla u)\|^2_2>\frac{p\gamma_pNc^2}{2(p\gamma_p-2)}=k_0.
\end{equation*}
In fact, since $c<c_0$, we have
\begin{equation*}
\begin{aligned}
\frac{p\gamma_p-2}{p\gamma_p}e^{2\bar{t}}\|\nabla u\|^2_2-\frac{N}{2}c^2&=
\Big(\frac{2\|\nabla u\|^2_2}{p\gamma^2_p\mu\|u\|^p_p}\Big)^{\frac{2}{p\gamma_p-2}}\frac{p\gamma_p-2}{p\gamma_p}
\|\nabla u\|^2_2-\frac{N}{2}c^2\\
&\ge \Big(\frac{2}{p\gamma^2_p\mu} \Big)^{\frac{2}{p\gamma_p-2}}C^{-\frac{2p}{p\gamma_p-2}}(N,p)\frac{p\gamma_p-2}{p\gamma_p}
c^{-\frac{2p(1-\gamma_p)}{p\gamma_p-2}}-\frac{N}{2}c^2\\
&>0.\\
\end{aligned}
\end{equation*}
This gives $\|\bar{t}\star (\nabla u)\|^2_2>k_0$.
It follows that there exists $\delta>0$ such that for any $(\bar{t}-\delta,\bar{t})$
\begin{equation*}
\|t\star (\nabla u)\|^2_2>k_0.
\end{equation*}
From \eqref{z8}, we get that for any $(\bar{t}-\delta,\bar{t})$, $\Psi'_{u}(t)>0$ and thus $\Psi_{u}(t)$ increasing in $(\bar{t}-\delta,\bar{t})$.
\vskip1mm

Notice that
\[
\lim\limits_{s\to-\infty}\Psi_{u}(s)=+\infty\ \text{and}\  \lim\limits_{s\to+\infty}\Psi_{u}(s)=-\infty,
\]
we thus conclude that
$\Psi_{u}(s)$ has a local minimum point $s_{u}<\bar{t}$ and has a local maximum point $t_{u}>\bar{t}$. Since $s_u<\bar{t}$, by \eqref{z7},
we deduce that
\begin{equation}\label{z9}
2\|s_u\star(\nabla u)\|^2_2>p\gamma^2_p\mu\|s_u\star u\|^p_p.
\end{equation}
In addition, combining \eqref{z9} and $\Psi'_{u}(s_u)=0$, we have
\begin{equation*}
\Psi''_{u}(s_{u})=2\|s_u\star(\nabla u)\|^2_2-p\gamma^2_p\mu\|s_u\star u\|^p_p>0.
\end{equation*}
This implies that $s_u$ is a strict minimum point for $\Psi'_{u}(t)$ and $s_{u}\star u\in \mathcal{P}^+_{c}$.
\vskip1mm

Next, we prove that $s_u$ is unique. By contradiction, we assume that $s^*_{u}$ is another local minimum point. We then claim that $s^*_{u}<\bar{t}$.
If not, we assume $s^*_{u}\ge \bar{t}$, then
\begin{equation*}
\Psi''_{u}(s^*_{u})=2\|s^*_u\star(\nabla u)\|^2_2-p\gamma^2_p\mu\|s^*_u\star u\|^p_p<0,
\end{equation*}
which is a contradiction. Therefore, $s^*_{u}<\bar{t}$ and $s^*_{u}\star u\in \mathcal{P}^+_{c}$. Moreover, there exists a another point,
$s^0_u<\bar{t}$ such that $s^0_u$ is a local maximum point for $\Psi_{u}$. From \eqref{z7}, we have $\Psi''_{u}(s^0_{u})>0$, which is a contradiction.
Thus $s_u$ is unique. Similar to the proof of $s_u$ is unique, we can prove $t_u$ is also unique.
\vskip2mm

$(ii)$ Finally, applying the Implicit Function Theorem to the $\mathcal{C}^1$ function $g(s,u) : \R \times S(c) \rightarrow \R$
defined by $g_{u}(s)=g(s,u)= \Psi_{u}'(s)$. Therefore, we have that
$u \mapsto s_{u}$ is of class $\mathcal{C}^1$ because $g_{u}(s_{u}) =0$ and $\partial_s g_{u}(s_{u}) = \Psi_{u}''(s_{u})>0$.
Similarly, we can prove that $u \mapsto t_{u} \in \mathbb{R}$ is of class $\mathcal{C}^1$.
\qed

\vskip1mm

\begin{lem}\label{lem4.5}
Let $\mu>0$ and $2+\frac{4}{N}<p<2^*$. Suppose that $0<c<c_0$. Then
$$m_c:=\inf_{u \in V_{k_0}} I(u)=m^+(c),$$
where $V_{k_0}$ is defined in \eqref{y1}.
\end{lem}
\noindent{\bf Proof.}
From the definition of $m^+(c)$, we can deduce that $\mathcal{P}^+_c\subset V_{k_0}$, then $m_c\le m^+(c)$.
On the other hand, if $u\in V_{k_0}$, $s_{u}\star u \!\in\!\mathcal{P}^{+}_{c}\subset V_{k_0}$, we have
\begin{equation*}
I(s_{u}\star u)=\min\big\{I(s\star u)~|~s\in \R~\text{and}~\|s\star(\nabla u)\|^2_2<k_0 \big\}\le I(u),
\end{equation*}
it follows that $m^+(c)=\inf\limits_{\mathcal{P}^+_c}I\le m_c.$
\qed

Now, we define
\[
\mathcal{P}^+_{r,c}:=\mathcal{P}^+_{c}\bigcap H^1_r(\R^N)\ \
\text{and}\ \
\mathcal{P}^-_{r,c}:=\mathcal{P}^-_{c}\bigcap H^1_r(\R^N).
\]

\begin{lem}\label{lem4.6}
Let $\mu>0$ and $2+\frac{4}{N}<p<2^*$. If $0<c< c_0$, then
\begin{equation*}
m^{+}(c):=\inf_{u\in\mathcal{P}^{+}_{c}} I(u)=\inf_{u\in \mathcal{P}^{+}_{r,c}} I(u)\ \ \text{and}\ \
m^{-}(c):=\inf_{u\in\mathcal{P}^{-}_{c}} I(u)=\inf_{u\in \mathcal{P}^{-}_{r,c}} I(u).
\end{equation*}
Moreover, $m^{+}(c)$, $m^{-}(c)$ are achieved by some Schwarz symmetric functions.
\end{lem}
\noindent{\bf Proof.}
Since $\mathcal{P}^{+}_{r,c}\subset \mathcal{P}^{+}_{c}$, $\mathcal{P}^{-}_{r,c}\subset \mathcal{P}^{-}_{c}$, we have
\[
\inf\limits_{u\in \mathcal{P}^{+}_{r,c}} I(u)\ge\inf\limits_{u\in\mathcal{P}^{+}_{c}} I(u)\ \ \text{and}\ \ \inf\limits_{u\in \mathcal{P}^{-}_{r,c}} I(u)\ge\inf\limits_{u\in\mathcal{P}^{-}_{c}} I(u).
\]
On the other hand, by Lemma \ref{lem4.3}, for each $u\in S(c)$, there exist $s_{u},t_u\in \R$ such that $s_{u}\star u\in\mathcal{P}^{+}_{c},t_{u}\star u\in\mathcal{P}^{-}_{c}$, and
\begin{equation*}
\inf_{u\in\mathcal{P}^{+}_{c}} I(u)=\inf_{u\in S(c)}\min _{-\infty<t\le s_{u}} I(t\star u),
\end{equation*}
\begin{equation*}
\inf_{u\in\mathcal{P}^{-}_{c}} I(u)=\inf_{u\in S(c)}\max _{s_{u}<t\le t_{u}} I(t\star u).
\end{equation*}
For $u\in S(c)$, let $u^*\in S_r(c)$ be the Schwarz rearrangement of $|u|$. By \cite[Chapter 3]{LL}, we have for all $t>0$,
$I(t\star u^*)\le I(t\star u)$. Recall that $\Psi'_{u}(t)=P(t \star u)$, we have
\begin{equation*}
\lim_{s\to -\infty}\Psi'_{u^*}(s)\le \lim_{s\to -\infty}\Psi'_{u}(s)=+\infty \quad \text{and} \quad \Psi''_{u^*}(t)\le \Psi''_{u}(t), \ \ \text{for each}\ \ t\in \R.
\end{equation*}
This implies that $-\infty<s_{u}\le s_{u^*}<t_{u^*}\le t_{u}$. Therefore, we get that
\begin{equation*}
\min _{-\infty<t<s_{u^*}} I(t\star u^*)\le \min _{-\infty<t<s_{u}} I(t\star u)
\end{equation*}
and
\begin{equation*}
\max _{s_{u^*}<t\le t_{u^*}} I(t\star u^*)\le \max _{s_{u}<t\le t_{u}} I(t\star u).
\end{equation*}
Therefore, $\inf\limits_{u\in \mathcal{P}^{+}_{r,c}} I(u)\le\inf\limits_{u\in\mathcal{P}^{+}_{c}} I(u)$
and $\inf\limits_{u\in \mathcal{P}^{-}_{r,c}} I(u)\le\inf\limits_{u\in\mathcal{P}^{-}_{c}} I(u)$.

Now if $u_0\in \mathcal{P}^+_{c}$ such that $I(u_0)=\inf\limits_{u\in \mathcal{P}^+_{c}}I(u)$, we see that $u^*_0$,
the Schwarz rearrangement of $|u_0|$, belongs to $\mathcal{P}^+_{r,c}$. In fact, if $\|\nabla u^*_0\|^2_2<\|\nabla u_0\|^2_2$,
then $I(u^*_0)<I(u_0)$. We get
\begin{equation*}
\begin{aligned}
\inf_{u\in\mathcal{P}^{+}_{c}} I(u)&=\inf_{u\in S(c)}\min _{-\infty<t\le s_{u}} I(t\star u)\\
&\le \min_{-\infty<t\le s_{u^*_0}} I(t\star u^*_0)< \min_{-\infty<t\le s_{u_0}} I(t\star u_0)\\
&=\inf_{u\in\mathcal{P}^{+}_{c}} I(u),
\end{aligned}
\end{equation*}
which is a contradiction. Therefore, $u^*_0\in \mathcal{P}^{+}_{r,c}$ and $I(u^*_0)=I(u_0)$.
Similarly, we can prove that there exists $u_0\in \mathcal{P}^-_{c}$ such that $I(u_0)=\inf\limits_{u\in \mathcal{P}^-_{c}}I(u)$.
\qed

\begin{lem}\label{lem4.7}
Let $\mu>0$ and $2+\frac{4}{N}<p<2^*$. If $0<c<c_0$, then $m^{\pm}(c)$ can be achieved by some $u^{\pm}_c$, which is positive, radially symmetric and decreasing in $r=|x|$.
\end{lem}
\noindent{\bf Proof.}
By Lemma \ref{lem4.2} and Lemma \ref{lem4.5}, the proof for $m^{+}(c)$ is completed. We only give the proof for $m^-(c)$.
Let $\{u_n\}\subset \mathcal{P}^{-}_{c}$ be a minimizing sequence of $m^-(c)$. From Lemma \ref{lem4.6}, then by taking $|u_{n}|$ and adapting the Schwarz symmetrization to
$|u_{n}|$ if necessary, we can obtain a new minimizing sequence (up to a subsequence),
such that $u_n$ are all nonnegative, radially symmetric and decreasing in $r = |x|$. Since $u_n\in \mathcal{P}^{-}_{c}$, we have
\begin{equation*}
\begin{aligned}
I(u_n)&= \Big(\frac{1}{2}-\frac{1}{p\gamma_p}\Big)\|\nabla u_n\|^2_2-\frac{1}{2}\int_{\R^N}u^2_n\log u^2_n\mathrm{d}x+\frac{p}{2(p-2)}c^2\\
&\ge \Big(\frac{1}{2}-\frac{1}{p\gamma_p}\Big)\|\nabla u_n\|^2_2 +\frac{1}{2}\int_{\R^N}A(|u_n|)\mathrm{d}x-\frac{1}{2}\int_{\R^N}B(|u_n|)\mathrm{d}x\\
&\ge \Big(\frac{1}{2}-\frac{1}{p\gamma_p}\Big)\|\nabla u_n\|^2_2 +\frac{1}{2}\int_{\R^N}A(|u_n|)\mathrm{d}x-\frac{K_q}{2}\int_{\R^N}|u_n|^q\mathrm{d}x,
\end{aligned}
\end{equation*}
where $2<q<2+\frac{4}{N}$.
Thus, by the Gagliardo-Nirenberg inequality \eqref{b2}, $\{u_n\}$ is bounded in $H^1(\R^N)$ and $\int_{\R^N}A(|u_n|)\mathrm{d}x$ is bounded. By Lemma \ref{lem2.1},
we get $\{u_n\}$ is bounded in $W_r$. By Lemma \ref{lem2.2}, $u_n\rightharpoonup u_0$ weakly in $W_r$, $u_n\to u_0$ strongly in $L^2(\R^N)$ and $u_n\to u_0$ a.e. in $\R^N$,
and thus $u_0\in S(c)$, $P(u_0)\le 0$.
From Lemma \ref{lem4.6}, there exists $s_0\in \R$ such that $s_0\star u_0\in \mathcal{P}^{-}_c$,
\begin{equation*}
m^{-}(c)+o_n(1)=I(u_n)\ge I(s_0\star u_n)\ge I(s_0\star u_0)+o_n(1)\ge m^{-}(c).
\end{equation*}
$I(s_0\star u_0)=m^{-}(c)$, $A(|u_n|)\to A(|s_0\star u_0|)$ in $L^1(\R^N)$. It thus follows from Lemma \ref{lem2.1}
that $u_n\to s_0\star u_0$ in $V$ as $n\to \infty$. Therefore, $u_n\to s_0\star u_0$ in $W$.
Finally, it follows from the proof in Lemma \ref{lem2.3} that $ s_0\star u_0> 0$.
We deduce that $m^{-}(c)$ is attained by $s_0\star u_0$ which is positive radially symmetric and decreasing in $r=|x|$.
\qed

\vskip 0.12truein

{\bf Proof of Theorem \ref{th1.2}.}
The proof follows directly from Lemma \ref{lem4.7}.

\section{Sobolev-critical case}

In this section, we study the case $\alpha=1$, $\mu>0$ and $p=2^*$.
Denote
\[
\mathcal{S}=\inf _{u \in {D}^{1,2}(\mathbb{R}^{N})\setminus \{0\} }    \frac{\left\| \nabla u\right\|_2^{2}}{\|u\|_{2^*}^{2}}.
\]
Then the following Sobolev inequality hold
\begin{equation}\label{Sobolev-ineq}
\mathcal{S}\|u\|_{2^*}^2\leq \|\nabla u\|_2^2,\quad \text{for each}\ \ u\in D^{1,2}(\R^N),
\end{equation}
where $D^{1,2}(\R^N)$ is the completion of $\mathcal{C}_c^{\infty}(\R^N)$ with respect to the norm $\|u\|_{D^{1,2}}:=\|\nabla u\|_2$.
For $c>0$, $N\ge 3$, we denote
\begin{equation}\label{z11}
k_0=\frac{N^2c^2}{4},  \qquad  c_0=\left[\frac{N^2-2N}{4\mu}\Big(\frac{4\mathcal{S}}{N^2}\Big)^{\frac{2^*}{2}}\right]^{\frac{1}{2^*-2}}.
\end{equation}

\subsection{Existence of a local minima on $S(c)$}
\begin{lem}\label{lem5.1}
Suppose that $\mu>0$, $c>0$ and $p=2^*$. If $P(u)\le 0$ and $\|\nabla u\|^2_2=k_0$, then $c\ge c_0$. Moreover, if $P(u)\le 0$ and $c<c_0$, then $\|\nabla u\|^2_2\neq k_0$.
\end{lem}
\noindent{\bf Proof.}
Since $P(u)\le 0$, we have
\begin{equation*}
\|\nabla u\|^2_2\le \mu \|u\|^{2^*}_{2^*}+\frac{N}{2}c^2.
\end{equation*}
By the Sobolev inequality \eqref{Sobolev-ineq}, it is easy to verify that
\begin{equation*}
\|\nabla u\|^2_2\le \mu  \mathcal{S}^{-\frac{2^*}{2}}\|\nabla u\|^{2^*}_2+\frac{N}{2}c^2.
\end{equation*}
By using the fact $\|\nabla u\|^2_2=k_0$ that
\begin{equation*}
\frac{N^2-2N}{4}\le \mu  \mathcal{S}^{-\frac{2^*}{2}}\Big(\frac{N^2}{4}\Big)^{\frac{2^*}{2}} c^{2^*-2},
\end{equation*}
it follows that $c\ge c_0$.
Therefore, we deduce that if $P(u)\le 0$ and $c<c_0$, then $\|\nabla u\|^2_2\neq k_0$.
\qed

\vskip 0.12truein

Define
\begin{equation}\label{v7}
V_{k_0}:=\left\{u\in S(c)~|~ \|\nabla u\|^2_2<k_0\right \}.
\end{equation}
where $k_0$ is defined by \eqref{z11}. Thus, for any $0<c\le c_0$,
we consider the following minimization problem
\begin{equation}\label{z21}
m_c := \inf_{u \in V_{k_0}} I(u).
\end{equation}

\begin{lem}\label{lem5.2}
Let $\mu>0$ and $p=2^*$, $N\ge 3$. Suppose that $0<c\le c_0$.
Then $m_c$ is achieved by some $u$, which
is a positive, radially symmetric critical point of $I$ on $S(c)$.
\end{lem}
\noindent{\bf Proof.}
Let $\{w_n\}$ be a minimizing sequence for $m_c$. Using the symmetric decreasing rearrangement, we obtain a minimizing sequence $\{w_n^*\}\subset H^1_r(\R^N)\cap V_{k_0}$, $w_n^*$ is radial decreasing for each $n$. Thus, by the Ekeland's variational principle,
we can obtain a nonnegative radial Palais-Smale sequence $\{u_n\}$ for $I|_{S(c)}$ at level $m_c$, i.e.
\[
\lim\limits_{n\to \infty} I(u_n)=m_c\ \ \text{and}\ \
dI\vert_{S(c)}(u_n)\to 0\ \ \text{as}\ \ n\to \infty.
\]
Since
\begin{equation*}
\begin{aligned}
I(u_n)&= \frac{1}{2}\|\nabla u_n\|^2_2-\frac{1}{2}\int_{\R^N}u^2_n\log u^2_n\mathrm{d}x-\frac{\mu}{2^*}\int_{\R^N}|u_n|^{2^*}\mathrm{d}x+\frac{1}{2}c^2\\
&\ge \frac{1}{2}\|\nabla u_n\|^2_2 +\frac{1}{2}\int_{\R^N}A(|u_n|)\mathrm{d}x-\frac{1}{2}\int_{\R^N}B(|u_n|)\mathrm{d}x-\frac{\mu}{2^*}\mathcal{S}^{-\frac{2^*}{2}}k^{\frac{2^*}{2}}_0,\\
\end{aligned}
\end{equation*}
by the same method as the proof of Lemma \ref{lem4.7}, we can show that $\{u_n\}$ is bounded in $W_r$, then $u_n\rightharpoonup u$ weakly in $W_r$, $u_n\to u$ strongly in $L^2(\R^N)$.
In fact, by the Lagrange multiplier's rule (see \cite[Lemma 3]{BL}) there is $\{\lambda_{n}\}\subset \R$ such that
\begin{equation}\label{b14}
\int_{\R^N}\big( \nabla u_{n}\nabla\phi+\lambda_{n}u_{n}\phi-\mu|u_{n}|^{p-2}u_{n}\phi- u_n\phi\log u^2_n\big)\mathrm{d}x=o_n(1)\|\phi\|_{W},
\end{equation}
for each $\phi\in W$. In particular, if we take $\phi=u_{n}$, we can prove that $\{\lambda_{n}\}$ is bounded.
Hence, up to a subsequence, we assume $\lambda_{n}\to \lambda \in \R$. Passing to the limit in \eqref{b14}, we get that  $u$ satisfies
\begin{equation*}
-\Delta u+\lambda u=u\log u^2+\mu|u|^{2^*-2}u.
\end{equation*}
Set $v_n := u_n-u$, we get
\begin{equation*}
\int_{\R^N}|\nabla v_n|^2+A(|v_n|)\mathrm{d}x-\mu\int_{\R^N}|v_n|^{2^*}\mathrm{d}x=o_n(1).
\end{equation*}
Therefore
\begin{equation*}
\begin{aligned}
m_c&\le I(u)\le\frac{1}{2}\int_{\R^N}|\nabla u|^2+u^2\mathrm{d}x-\frac{1}{2}\int_{\R^N}u^2\log u^2\mathrm{d}x-\frac{\mu}{2^*}\int_{\R^N}|u|^{2^*}\mathrm{d}x\\
&\qquad \qquad +\frac{1}{2}\int_{\R^N}|\nabla v_n|^2+A(|v_n|)\mathrm{d}x-\frac{\mu}{2^*}\int_{\R^N}|v_n|^{2^*}\mathrm{d}x\\
&= \frac{1}{2}\int_{\R^N}|\nabla u_n|^2+u^2_n\mathrm{d}x-\frac{1}{2}\int_{\R^N}u^2_n\log u^2_n\mathrm{d}x-\frac{\mu}{2^*}\int_{\R^N}|u_n|^{2^*}\mathrm{d}x+o_n(1)\\
&= \liminf_{n \to \infty}I(u_n)=m_c.
\end{aligned}
\end{equation*}
Next, we only need to prove $u\in V_{k_0}$. In fact, if $\|\nabla u\|^2_2=k_0$ and $0<c<c_0$, it follows from Lemma \ref{lem5.1} that $P(u)>0$,
then there exists $t_0<0$ such that $t_0\star u\in V_{k_0}$ and
$I(t_0\star u)<I(u)=m_c.$ This is a contradiction. On the other hand, if $\|\nabla u\|^2_2=k_0$ and $c=c_0$, the proof is similar to Lemma \ref{lem4.2}.
Therefore, $u\in V_{k_0}$ and $I(u)=m_c$.
According to the proof in Lemma \ref{lem3.1}, we conclude that $u>0$.
\qed

\subsection{Multiplicity of solutions}

In this subsection, we are interested in the multiplicity of solutions. The following two Lemmas are similar to Lemma \ref{lem4.3} and Lemma \ref{lem4.4}, so we omit the proof.

\begin{lem}\label{lem5.3}
Let $\mu>0$ and $p=2^*$. If $0<c<c_0$, then $\mathcal{P}^0_{c}=\emptyset$, and the set $\mathcal{P}_{c}$ is a $\mathcal{C}^1$-submanifold of codimension 1 in $S(c)$.
\end{lem}

\begin{lem}\label{lem5.4} If $p=2^*$ and $0<c<c_0$, for $u\in S(c)$, then the function
$\Psi_{u}(s)$ has exactly two critical points $s_{u}\!<\!t_{u}\in\R$. Moreover,

\begin{itemize}
\item  [$(i)$] $s_{u}\star u \!\in\!\mathcal{P}^{+}_{c}$ and $s_u$ is a strict local minimum point for $\Psi_{u}$. $t_{u}\star u \!\in\!\mathcal{P}^{-}_{c}$ and $t_u$ is a strict local maximum point for $\Psi_{u}$.

\item  [$(ii)$] The maps $u \mapsto s_{u} \in \mathbb{R}$
and $u \mapsto t_{u} \in \mathbb{R}$ are of class $\mathcal{C}^1$.
\end{itemize}
\end{lem}

In view of Lemma \ref{lem5.4}, we can define
\begin{equation*}
m^+(c):=\inf_{u\in\mathcal{P}^+_c}I(u)\quad \text{and} \quad m^-(c):=\inf_{u\in\mathcal{P}^-_c}I(u).
\end{equation*}

Based on $\sigma$-homotopy stable family of compact subsets of $\mathcal{P}_c$, we aim to construct bounded
Palais-Smale sequences on a manifold by the Ghoussoub minimax principle \cite[Theorem 3.2]{GN}, (see also \cite[Section2.2]{JJ}, \cite[Lemma 3.7]{JL} and \cite[Section 5]{SN}).

\vskip1mm
Define the functions
\begin{equation*}
\begin{aligned}
J^{+}: S(c)\to \R, \quad J^+:=I(s_u\star u),\\
J^{-}: S(c)\to \R, \quad J^-:=I(t_u\star u).\\
\end{aligned}
\end{equation*}
Similar as Lemma \ref{lem5.4} $(ii)$, we can prove the maps $u\mapsto s_{u}$ and $u\mapsto t_{u}$ are of class $\mathcal{C}^1$,
thus the functionals $J^+$, $J^-$ are of class $\mathcal{C}^1$.
\begin{lem}\label{lem5.5}
For any $0<c<c_0$, we have that
\begin{equation*}
dJ^+[\psi]=dI(s_{u}\star u)[s_{u}\star \psi]\quad \text{and} \quad dJ^{-}[\psi]=dI(t_{u}\star u)[t_{u}\star \psi],
\end{equation*}
for any $u\in S(c)$ and $\psi\in T_{u}S(c)$.
\end{lem}
\noindent{\bf Proof.}
We first give the proof for $J^+$. Let $\psi\in S(c)$, and $\psi=h'(0)$ where $h:(-\epsilon,\epsilon)\mapsto S(c)$ is a
$\mathcal{C}^1$-curve with $h(0)=u$. We consider the incremental quotient
\begin{equation*}
\frac{J^+(h(t))-J^+(h(0))}{t}=\frac{I(s_{t}\star h(t))-I(s_{0}\star h(0))}{t},
\end{equation*}
where $s_{t}:=s_{h(t)}$ and $s_{0}:=s_{h(0)}=s_{u}$. It follows from Lemma \ref{lem5.4} that $s_0$ is a strict local minimum of
$s\mapsto I(s\star u)$, 
we have
\begin{equation*}
\begin{aligned}
&I(s_{t}\star h(t))-I(s_{0}\star h(0))\ge I(s_{t}\star h(t))-I(s_{t}\star h(0))\\
&=\frac{e^{2s_t}}{2}\big(\|\nabla h(t)\|^2_2-\|\nabla h(0)\|^2_2\big)+\big(\frac{1}{2}-\frac{Ns_{t}}{2}\big)\big(\|h(t)\|^2_2-\|h(0)\|^2_2\big)\\
&\quad-\frac{1}{2}\big(\int_{\R^N}|h(t)|^2\log|h(t)|^2\mathrm{d}x-\int_{\R^N}|h(0)|^2\log|h(0)|^2\mathrm{d}x\big)
-\frac{\mu e^{2^*s_{t}}}{2^*}\big(\|h(t)\|^{2^*}_{2^*}-\|h(0)\|^{2^*}_{2^*}\big)\\
&=e^{2s_t}\int_{\R^N}\nabla h(\tau_1 t)\nabla h'(\tau_1t)\mathrm{d}x+\big(1-Ns_{t}\big)\int_{\R^N}h(\tau_2t)h'(\tau_2t)\mathrm{d}x\\
&\quad-\int_{\R^N}h(\tau_3t)h'(\tau_3t)\log|h(t)|^2+|h(0)|^2\frac{h'(\tau_4t)}{h(\tau_4t)}\mathrm{d}x
-\mu e^{2^*s_{t}}\int_{\R^N}|h(\tau_5t)|^{2^*-2}h(\tau_5t)h'(\tau_5t)\mathrm{d}x,\\
\end{aligned}
\end{equation*}
for some $\tau_1,\cdots,\tau_5\in(0,1)$.
In the same manner, we get
\begin{equation*}
\begin{aligned}
&I(s_{t}\star h(t))-I(s_{0}\star h(0))\le I(s_{0}\star h(t))-I(s_{0}\star h(0))\\
&\quad=e^{2s_0}\int_{\R^N}\nabla h(\tau_6 t)\nabla h'(\tau_6t)\mathrm{d}x+\big(1-Ns_{0}\big)\int_{\R^N}h(\tau_7t)h'(\tau_7t)\mathrm{d}x\\
&\quad\quad-\int_{\R^N}h(\tau_8t)h'(\tau_8t)\log|h(t)|^2+|h(0)|^2\frac{h'(\tau_9t)}{h(\tau_9t)}\mathrm{d}x\\
&\quad\quad-\mu e^{2^*s_{t}}\int_{\R^N}|h(\tau_{10}t)|^{2^*-2}h(\tau_{10}t)h'(\tau_{10}t)\mathrm{d}x,\\
\end{aligned}
\end{equation*}
for some $\tau_6,\cdots,\tau_{10}\in(0,1)$. Finally, we deduce that
\begin{equation*}
\begin{aligned}
\lim_{t\to 0}\frac{J^+(h(t))-J^+(h(0))}{t}&=\int_{\R^N}\nabla (s_{u}\star u)\nabla (s_{u}\star \psi)\mathrm{d}x
-\int_{\R^N}(s_{u}\star u)(s_{u}\star \psi)\log|s_{u}\star u|^2\mathrm{d}x\\
&\quad-\mu \int_{\R^N}|s_{u}\star u|^{2^*-2}(s_{u}\star u)(s_{u}\star \psi)=dI(s_{u}\star u)[s_{u}\star \psi],
\end{aligned}
\end{equation*}
for any $u\in S(c)$ and $\psi\in T_{u}S(c)$. The proof for $J^-$ is similar.
\qed

\vskip1mm
Let $\mathcal{G}$ be the set of all singletons belongs to $S(c)$ and boundary $B=\emptyset$. It is a homotopy stable family of compact
subset of $S(c)$ (without boundary) in the sense of \cite[Definition 3.1]{GN}.
By Lemma \ref{lem4.6}, we have
\begin{equation*}
e^{+}_{\mathcal{G}}:=\inf_{A\in\mathcal{G}}\max_{u\in A}J^{+}(u)=\inf_{u\in S(c)\cap H^1_r(\R^N)}J^{+}(u)
=\inf_{u\in \mathcal{P}^{+}_{r,c}}I(u)=\inf_{u\in \mathcal{P}^{+}_{c}}I(u),
\end{equation*}
\begin{equation*}
e^{-}_{\mathcal{G}}:=\inf_{A\in\mathcal{G}}\max_{u\in A}J^{-}(u)=\inf_{u\in S(c)\cap H^1_r(\R^N)}J^{-}(u)
=\inf_{u\in \mathcal{P}^{-}_{r,c}}I(u)=\inf_{u\in \mathcal{P}^{-}_{c}}I(u).
\end{equation*}

\begin{lem}\label{lem5.6}
For any $0<c<c_0$, there exists a Palais-Smale sequence $\{u_n\}^{+}\subset \mathcal{P}^{+}_{c}$ (or $\{u_n\}^{-}\subset \mathcal{P}^{-}_{c}$) for $I$ restricted to
$S(c)\cap H^1_r(\R^N)$ at level $e^{+}_{\mathcal{G}}$ (or $e^{-}_{\mathcal{G}}$, respectively).
\end{lem}
\noindent{\bf Proof.}
We first prove the case of $e^{+}_{\mathcal{G}}$, the one for $e^{-}_{\mathcal{G}}$ is almost identical. By the definition of $e^{+}_{\mathcal{G}}$, there exists $\{E_n\}\subset \mathcal{G}$
such that
\begin{equation*}
\max_{u\in E_n}J^+(u)<e^{+}_{\mathcal{G}}+\frac{1}{n}.
\end{equation*}
Define
\begin{equation*}
\eta:[0,1]\times S(c)\mapsto S(c), \quad \eta(t,u)=(1-t+ts_{u})\star u.
\end{equation*}
By using the definition of $\mathcal{G}$, we have
\begin{equation*}
\bar{E}_n:=\eta(\{1\}\times E_n)=\{s_{u}\star u: u\in E_n\}\in \mathcal{G}.
\end{equation*}
It follows from Lemma \ref{lem5.4} that $\bar{E}_n\subset \mathcal{P}^+_c$ for all $n\in \N$. Let $v\in \bar{E}_n$, i.e. $v=s_{u}\star u$
for some $u\in E_n$, then $J^+(v)=J^+(u)$, and finally that
\begin{equation*}
\max_{v\in \bar{E}_n}J^+(v)=\max_{u\in E_n}J^+(u).
\end{equation*}
Using the terminology in \cite[Definition 3.1]{CT},
it means that $\bar{E}_n$ is a homotopy stable family of compact subset of $S(c)\cap H^1_r(\R^N)$.
Therefore, $\{\bar{E}_n\}$ is another sequence sets such that $\lim\limits_{n}\sup\limits_{\bar{E}_n}J^+=e^+_{\mathcal{G}}$.
Denoting by $\|\cdot\|_{\ast}$ the dual norm of $(T_{v_n}S(c))^{*}$,
we can apply \cite[Theorem 3.2]{GN} with
the minimizing sequence $w_n$ for $J^+$ on $S(c)$ at level $e^+_{\mathcal{G}}$ such that
\begin{equation}\label{v1}
(i)~\lim_{n\to \infty}J^+(w_n)=e^+_{\mathcal{G}},\ \ (ii)~\lim_{n\to \infty}\|dJ^+(w_n)\|_{*}=0,\ \ (iii)~{dist}_{W_r}(w_n, \bar{E}_n)\to 0\ \ \text{as}\ n\to \infty.\\
\end{equation}
Denote $v_n:=s_{w_n}\star w_n\in \mathcal{P}^+_{c}$, we claim that $\{s_{w_n}\}$ is bounded. Notice that $$e^{2s_{w_n}}=\frac{\|\nabla v_n\|^2_2}{\|\nabla w_n\|^2_2}.$$
In fact, $I(v_n)=J^+(w_n)\to e^+_{\mathcal{G}}=m^+(c)$, we get
from that $v_n$ is uniformly bounded in $W_r$. We deduce that $\{\bar{E}_n\}$ is uniformly bounded in $W_r$, and thus from
$dist_{W_r}(w_n, \bar{E}_n)\to 0$ as $n\to \infty$, it gives that $\sup_{n}\|\nabla w_n\|^2_2<\infty$. Since $\bar{E}_n$ is compact for every
$n\in\N$, there exists $\tilde{v}_n\in \bar{E}_n$ such that $dist_{W_r}(w_n, \bar{E}_n)=\|\tilde{v}_n-w_n\|_{W_r}$, and then
$\|\nabla w_n\|^2_2\ge\|\nabla\tilde{v}_n\|^2_2-\|\nabla(\tilde{v}_n-w_n)\|^2_2\ge \frac{\delta}{2}$ for a $\delta>0$. Therefore, $\{s_{w_n}\}$ is bounded.
\vskip1mm

Next, we show that $\{v_n\}\subset \mathcal{P}^+_c$ is a Palais-Smale sequence for $I$ restricted on $S(c)$ at level $e^+_{\mathcal{G}}$.
From Lemma \ref{lem2.3}, we get that $T_{w_n}S(c)\to T_{v_n}S(c)$ defined $\psi\to s_{w_n}\star \psi$ is an isomorphism.
By \eqref{v1}, we get
\begin{equation*}
\|dI|_{S(c)}(v_n)\|_{*}=\sup_{\|\psi\|\le 1,\psi\in T_{u}S(c)}|dI(v_n)[\psi]|
=\sup_{\|\psi\|\le 1,\psi\in T_{u}S(c)}|dJ^{+}(w_n)[(-s_{w_n})\star\psi]|,
\end{equation*}
since $s_{w_n}$ is bounded and $\|(-s_{w_n})\star\psi\|_{W_r}\le C$.
We conclude that $\{v_n\}\subset \mathcal{P}^+_c$ is a Palais-Smale sequence for $I$ restricted to $S(c)$ at level $e^+_{\mathcal{G}}$.
\qed

\vskip 0.12truein

Applying a similar arguments as the proof of Lemma \ref{lem4.5}. We can deduce that
\begin{lem}\label{lem5.7}
Let $\mu>0$ and $p=2^*$. Suppose that $0<c<c_0$. Then
\begin{equation}\label{eqs5.7-7}
m_c:=\inf_{u \in V_{k_0}} I(u)=m^+(c),
\end{equation}
where $V_{k_0}$ is defined in \eqref{v7}.
\end{lem}

\begin{lem}\label{Lem5.8}
Let $\mu>0$, and $p=2^*$. Suppose that $0<c<c_0$. Then $m^{+}(c)$ is achieved.
\end{lem}
\noindent{\bf Proof.}
By Lemma \ref{lem5.7}, this follows by the same method as in Lemma \ref{lem5.2}.
Applying Lemma \ref{lem5.6}, we deduce that there exists a Palais-Smale sequence $\{u_n\}\subset \mathcal{P}^+_{c}$ for $I$ restricted
to $S(c)$ at level $e^+_{\mathcal{G}}=m^+(c)$.
Obviously,
\begin{equation*}
\begin{aligned}
I(u_n)&= \frac{1}{2}\|\nabla u_n\|^2_2-\frac{1}{2}\int_{\R^N}u^2_n\log u^2_n\mathrm{d}x
-\frac{\mu}{2^*}\int_{\R^N}|u_n|^{2^*}\mathrm{d}x+\frac{1}{2}c^2\\
&\ge \Big(\frac{1}{2}-\frac{1}{2^*}\Big)\|\nabla u_n\|^2_2 +\frac{1}{2}\int_{\R^N}A(|u_n|)-B(|u_n|)\mathrm{d}x+\frac{N}{4}c^2.\\
\end{aligned}
\end{equation*}
The proof is similar to that of Lemma \ref{lem4.7}, we get $\{u_n\}$ is bounded in $W_r$, $u_n\rightharpoonup u$ weakly in $W_r$, $u_n\to u$ strongly in $L^2(\R^N)$.
In fact, by the Lagrange multiplier's rule (see \cite[Lemma 3]{BL}) there is $\{\lambda_{n}\}\subset \R$ such that
\begin{equation}\label{b4-1}
\int_{\R^N} \big(\nabla u_{n}\nabla\phi+\lambda_{n}u_{n}\phi-\mu|u_{n}|^{2^*-2}u_{n}\phi- u_n\phi\log u^2_n\big)\mathrm{d}x=o_n(1)\|\phi\|_{W},\ \ \text{for each}\ \ \phi\in W.
\end{equation}
In particular, if we take $\phi=u_{n}$, we then conclude $\lambda_{n}$ is bounded. Hence, up to a subsequence, $\lambda_{n}\to \lambda \in \R$.
Passing to the limit in \eqref{b4-1}, we deduce that $u$ satisfies
\begin{equation*}
-\Delta u+\lambda u=u\log u^2+\mu|u|^{2^*-2}u,\ \ x\in \R^N.
\end{equation*}
It follows that $P(u)=0$.
Then
\begin{equation*}
\begin{aligned}
m^+(c)&\le I(u)\le\Big(\frac{1}{2}-\frac{1}{2^*}\Big)\int_{\R^N}|\nabla u|^2\mathrm{d}x+\frac{N}{4}c^2-\frac{1}{2}\int_{\R^N}u^2\log u^2\mathrm{d}x\\
&\le \Big(\frac{1}{2}-\frac{1}{2^*}\Big)\liminf_{n\to \infty}\int_{\R^N}|\nabla u_n|^2\mathrm{d}x+\frac{1}{2}\liminf_{n\to \infty}\int_{\R^N}A(|u_n|)\mathrm{d}x\\
&\quad-\lim_{n\to \infty}\int_{\R^N}B(|u_n|)\mathrm{d}x+\frac{N}{4}c^2\\
&\le \liminf_{n \to \infty}I(u_n)=m^{+}(c).
\end{aligned}
\end{equation*}
Therefore, $I(u)=m^+(c)$.
\qed

\begin{lem}\label{lem5.9}
Let $p=2^*$, $\mu>0$. Then for $0<c<c_0$, we have
\begin{equation*}
m^-(c)<m^+(c)+\frac{\mu^{-\frac{N-2}{2}}}{N}\mathcal{S}^{\frac{N}{2}}.
\end{equation*}
\end{lem}
\noindent{\bf Proof.}
In fact, we assume that $u\in S(c)$, $I(u)=m^+(c)$. Let
\[
U_\varepsilon(x)=\frac{[N(N-2)\varepsilon^2]^{\frac{N-2}{4}}}{[\varepsilon^2+|x|^2]^{\frac{N-2}{2}}}
\]
be the positive solution of
\[
-\Delta v =|v|^{2^*-2}v,\ \ v\in D^{1,2}(\R^N).
\]
Let $\overline{U}_{\varepsilon}=\chi(x)U_{\varepsilon}$ where $\chi$ is a cut-off function such that $\chi(x)=1$ for $|x|\le 1$ and $\chi(x)=0$ for $|x|>2$, we have
\begin{equation}\label{equ3.4}
\|\nabla \overline{U}_{\varepsilon}\|^2_2=\mathcal{S}^{\frac{N}{2}}+O(\varepsilon^{N-2})\quad  \text{and} \quad \|\overline{U}_{\varepsilon}\|^{2^*}_{2^*}=\mathcal{S}^{\frac{N}{2}}+O(\varepsilon^N).
\end{equation}
Define $W_{\varepsilon,t}=u(\cdot-ne_1)+t\overline{U}_{\varepsilon}$, where $e_1=(1,0,\cdots,0)$,
and $\overline{W}_{\varepsilon,t}=s^{\frac{N-2}{2}}W_{\varepsilon,t}(sx)$. Since $u$ is a radial, non-increasing function, we know from \cite[Radial
Lemma A.IV]{BL1} that
\begin{equation*}
|u(z)|\le C_N|z|^{-\frac{N}{2}}c, \quad \text{for}\ ~ |z|\ge 1.
\end{equation*}
Thus, we obtain
\begin{equation*}
-\Delta u+q(x)u=0,\qquad u\in H^1(\R^N),
\end{equation*}
where $q(x)=\lambda-\log u^2-\mu|u|^{2^*-2}$. Observe that
$q(x)\ge \frac{1}{2}$ for $|x|>M$ if $M$ is large enough.
By \cite[Lemma 3.11]{BS}, there exists $a\in (0,\frac{1}{2})$ such that
\begin{equation*}
u(x)\le C e^{-\sqrt{1+a|x|^2}},\qquad \text{for~each}~|x|>M.
\end{equation*}
By direct calculations, we have
\begin{equation*}
\|\nabla \overline{W}_{\varepsilon,t}\|^2_2=\|\nabla W_{\varepsilon,t}\|^2_2, \qquad \|\overline{W}_{\varepsilon,t}\|^{2^*}_{2^*}=\| W_{\varepsilon,t}\|^{2^*}_{2^*},
\end{equation*}
and
$$\|\overline{W}_{\varepsilon,t}\|^2_2=s^{-2}\|W_{\varepsilon,t}\|^2_2.$$
We choose $s=\frac{\|W_{\varepsilon,t}\|_2}{c}$ such that $\overline{W}_{\varepsilon,t}\in S(c)$. By Lemma \ref{lem4.3},
there exists $\tau_{\varepsilon,t}\in \R$ such that $\tau_{\varepsilon,t}\star \overline{W}_{\varepsilon,t}\in \mathcal{P}^{-}_{c}$. Thus
\begin{equation}\label{equ3.5}
e^{2\tau_{\varepsilon,t}}\|\nabla \overline{W}_{\varepsilon,t}\|^2_2=\mu e^{2^*\tau_{\varepsilon,t}}\|\overline{W}_{\varepsilon,t}\|^{2^*}_{2^*}
+\frac{N}{2}c^2.
\end{equation}
Since $u\in \mathcal{P}^+_{c}$, by Lemma \ref{lem4.3}, we deduce $\tau_{\varepsilon,0}>0$. By using \eqref{equ3.4} and \eqref{equ3.5},
we have $\tau_{\varepsilon,t}\to -\infty$ as $t\to +\infty$ and $\varepsilon>0$ small enough. It follows from Lemma \ref{lem4.3}
that $\tau_{\varepsilon,t}$ is unique. Moreover, $\tau_{\varepsilon,t}$ is continuous for $t$, 
and then there exists $t_{\varepsilon}$ such that $\tau_{\varepsilon,t_{\varepsilon}}=0$.
Therefore, we get
\begin{equation*}
 m^-(c)\le \sup_{t\ge 0}I\big(\overline{W}_{\varepsilon,t}\big).
\end{equation*}
There exists $t_0>0$ such that
\begin{equation*}
\begin{aligned}
I\big(\overline{W}_{\varepsilon,t}\big)&<m^+(c)+\frac{1}{N}\mu^{-\frac{N-2}{2}}\mathcal{S}^{\frac{N}{2}}-\sigma,
\end{aligned}
\end{equation*}
for $0<t<\frac{1}{t_0}$ and $t>t_0$ with $\sigma>0$.
\vskip1mm

Since the function $\overline{U}_\varepsilon$ is compacted supported in $B_2$, we have that, for $n$ large enough,
\begin{equation*}
\int_{\R^N}u(x-ne_1)\overline{U}_{\varepsilon}\mathrm{d}x=O(\varepsilon^{\frac{N-2}{2}}),
\end{equation*}
and then
\begin{equation}\label{equ3.6}
s^2=\frac{\|W_{\varepsilon,t}\|^2_2}{c^2}=1+\frac{2t}{c^2}\int_{\R^N}u(\cdot-ne_1)\overline{U}_{\varepsilon}\mathrm{d}x
+\frac{t^2\|\overline{U}_{\varepsilon}\|^2_2}{c^2},
\end{equation}
for $\frac{1}{t_0}\le t\le t_0$. Recall the definition of $A(s)$ in Section 2, for any $0\le a, b$ and $tb\ge \frac{1}{2}a$, $A(a+tb)\le A(a)+CtA(b)$ (see \cite[Lemma 1.3]{CT}). We deduce that
\begin{equation*}
\begin{aligned}
\int_{\R^N}|W_{\varepsilon,t}|^2\log|W_{\varepsilon,t}|^2\mathrm{d}x
&=\int_{\R^N}B(|u(x-ne_1)+t\overline{U}_{\varepsilon}|)\mathrm{d}x-\int_{\R^N}A(|u(x-ne_1)+t\overline{U}_{\varepsilon}|)\mathrm{d}x\\
&\ge \int_{\R^N}u^{2}\log u^{2}\mathrm{d}x +Ct\int_{\R^N}|\overline{U}_\varepsilon|^{2}\log|\overline{U}_\varepsilon|^{2}\mathrm{d}x.
\end{aligned}
\end{equation*}
By using \eqref{equ3.6}, we have
\begin{equation}\label{x1}
\begin{aligned}
I\big(\overline{W}_{\varepsilon,t}\big)&=\frac{1}{2}\int_{\R^N}|\nabla      W_{\varepsilon,t}|^2\mathrm{d}x+\frac{s^{-2}}{2}\int_{\R^N}|W_{\varepsilon,t}|^2\mathrm{d}x-\frac{s^{-2}\log|s^{N-2}|}{2}\int_{\R^N}|W_{\varepsilon,t}|^2\mathrm{d}x \\
&\quad-\frac{s^{-2}}{2}\int_{\R^N}|W_{\varepsilon,t}|^2\log|W_{\varepsilon,t}|^2\mathrm{d}x
-\frac{\mu}{2^*}\int_{\R^N}|W_{\varepsilon,t}|^{2^*}\mathrm{d}x  \\
&\le I(u)+\frac{t^2}{2}\int_{\R^N}|\nabla\overline{U}_\varepsilon|^2\mathrm{d}x
-\frac{t^{2^*}}{2^*}\mu\int_{\R^N}|\overline{U}_\varepsilon|^{2^*}\mathrm{d}x
+\frac{1}{2}\int_{\R^N}|\overline{U}_\varepsilon|^{2}\mathrm{d}x  \\
&\quad-Ct\int_{\R^N}|\overline{U}_\varepsilon|^{2}\log|\overline{U}_\varepsilon|^{2}\mathrm{d}x+\int_{\R^N}\nabla u(\cdot-ne_1)\nabla(tU_{\varepsilon})+u(\cdot-ne_1)tU_{\varepsilon}\mathrm{d}x\\
&\quad-\int_{\R^N}|u(\cdot-ne_1)|^{2^*-1}t\overline{U}_\varepsilon\mathrm{d}x -\frac{t}{c^2}\int_{\R^N}u(\cdot-ne_1)\overline{U}_\varepsilon\mathrm{d}x\|W_{\varepsilon,t}\|^2_2 \\
&\quad+\frac{t}{c^2}\int_{\R^N}u(\cdot-ne_1)\overline{U}_\varepsilon\mathrm{d}x\int_{\R^N}
|W_{\varepsilon,t}|^2\log|W_{\varepsilon,t}|^2\mathrm{d}x.\\
\end{aligned}
\end{equation}
Here, we used the fact that $f(s)=(1+s)^{2^*}\ge 1+s^{2^*}+2^*s^{2^*-1}\ge 0$ for all $s\ge 0$.
To complete the proof it suffices to show
\begin{equation}\label{x2}
\begin{aligned}
J_{\varepsilon,t}:&=\frac{1}{2}\int_{\R^N}|\overline{U}_\varepsilon|^{2}\mathrm{d}x-Ct\int_{\R^N}|\overline{U}_\varepsilon|^{2}\log|\overline{U}_\varepsilon|^{2}\mathrm{d}x
+\int_{\R^N}\nabla u(\cdot-ne_1)\nabla(tU_{\varepsilon})\mathrm{d}x\\
&\quad-\int_{\R^N}|u(\cdot-ne_1)|^{2^*-1}t\overline{U}_\varepsilon\mathrm{d}x +\int_{\R^N}u(\cdot-ne_1)tU_{\varepsilon}\mathrm{d}x\Big(1-\frac{1}{c^2}\|W_{\varepsilon,t}\|^2_2\Big)\\
&\quad+\frac{t}{c^2}\int_{\R^N}u(\cdot-ne_1)\overline{U}_\varepsilon\mathrm{d}x\int_{\R^N}|W_{\varepsilon,t}|^2\log|W_{\varepsilon,t}|^2\mathrm{d}x\\
&\le \frac{1}{2}\int_{\R^N}|\overline{U}_\varepsilon|^{2}\mathrm{d}x-Ct\int_{\R^N}|\overline{U}_\varepsilon|^{2}\log|\overline{U}_\varepsilon|^{2}\mathrm{d}x
+O(\varepsilon^{\frac{N}{2}})-O(\varepsilon^{\frac{N-2}{2}})+O(\varepsilon^{N-2})\\
&\quad+\frac{t}{c^2}\int_{\R^N}u(\cdot-ne_1)\overline{U}_\varepsilon\mathrm{d}x\int_{\R^N}|W_{\varepsilon,t}|^2\log|W_{\varepsilon,t}|^2\mathrm{d}x\\
&\le \frac{1}{2}\int_{\R^N}|\overline{U}_\varepsilon|^{2}\mathrm{d}x
-Ct\int_{\R^N}|\overline{U}_\varepsilon|^{2}\log|\overline{U}_\varepsilon|^{2}\mathrm{d}x+O(\varepsilon^{N-2})\\
&<0.\\
\end{aligned}
\end{equation}
Here, we used the fact that $\int_{\R^N}u(\cdot-ne_1)\overline{U}_\varepsilon\mathrm{d}x\le C\|\overline{U}_\varepsilon\|^2_2$,
because $u$ decays exponentially and $\overline{U}_\varepsilon$ has compact support.
\vskip1mm
If $N\ge 4$, we have, for $C_1,C_2>0$,
\begin{equation*}
\int_{\R^N}|\overline{U}_\varepsilon|^2\mathrm{d}x=
\begin{cases}
C_1\varepsilon^2|\ln\varepsilon|+O(\varepsilon^2), \qquad &\text{if}\ N=4,\\
C_2\varepsilon^2+O(\varepsilon^{N-2}), \qquad &\text{if}\ N\ge 5.\\
\end{cases}
\end{equation*}
From \cite[Lemma 3.4]{DPS}, we get
\begin{equation*}
\int_{\R^N}|\overline{U}_\varepsilon|^2\log|\overline{U}_\varepsilon|^2\mathrm{d}x\ge
\begin{cases}
C_3\varepsilon^2|\ln\varepsilon|^2+O(\varepsilon^2), \qquad &\text{if}\ N=4,\\
C_4\varepsilon^2|\ln\varepsilon|+O(\varepsilon^{N-2}), \qquad &\text{if}\ N\ge 5,\\
\end{cases}
\end{equation*}
where $C_3,C_4>0$. Therefore, $J_{\varepsilon,t}<0$. 

\vskip1mm

If $N=3$, by the definition of $\overline{U}_\varepsilon$, for $C_5>0$, we have
\begin{equation*}
\begin{aligned}
\int_{\R^3}|\overline{U}_\varepsilon|^2\log|\overline{U}_\varepsilon|^2\mathrm{d}x&=\sqrt{3}\int_{B_2}\chi^2\frac{\varepsilon}{\varepsilon^2+|x|^2}\
\log\big(\sqrt{3}\chi^2\frac{\varepsilon}{\varepsilon^2+|x|^2}\big)\mathrm{d}x\\
&=C_5 \varepsilon |\log \varepsilon|+O(\varepsilon).
\end{aligned}
\end{equation*}
Moreover, $\|\overline{U}_\varepsilon\|^2_2=O(\varepsilon)$. Therefore, $J_{\varepsilon,t}<0$.

\vskip1mm

Combining \eqref{x1} with \eqref{x2}, we have,
\begin{equation*}
I\big(\overline{W}_{\varepsilon,t}\big)< m^+(c)+\frac{1}{N}\mu^{-\frac{N-2}{2}}\mathcal{S}^{\frac{N}{2}},
\end{equation*}
for $\frac{1}{t_0}\le t_0\le t$.
Therefore, for $\varepsilon$ small enough,
we have
\begin{equation*}
I\big(\overline{W}_{\varepsilon,t}\big)<m^+(c)+\frac{1}{N}\mu^{-\frac{N-2}{2}}\mathcal{S}^{\frac{N}{2}}.
\end{equation*}
\qed

\begin{lem}\label{lem5.10}
Let $\mu>0$ and $p=2^*$. Suppose that $0<c<c_0$. Then $m^{-}(c)$ is achieved, there exists a second solution $u^-\in S(c)$ which satisfies
\begin{equation*}
I(u^{-})<m^+(c)+\frac{1}{N}\mu^{-\frac{N-2}{2}}\mathcal{S}^{\frac{N}{2}}.
\end{equation*}
\end{lem}
\noindent{\bf Proof.}
Similarly, by Lemma \ref{lem5.6}, we assume that there exists a Palais-Smale sequence $\{u_n\}\subset \mathcal{P}^-_{c}$ for $I$ restricted to $S(c)$ at level $m^-(c)$.
Moreover, $\{u_n\}$ is bounded in $W_r$, and $u_n\rightharpoonup u$ in $W_r$. We need to prove that $u_n\to u$ in $W_r$.
In fact, by the Lagrange multiplier's rule (see \cite[Lemma 3]{BL}) there is $\{\lambda_{n}\}\subset \R$ such that
\begin{equation}\label{b411}
\int_{\R^N} \big(\nabla u_{n}\nabla\phi+\lambda_{n}u_{n}\phi-\mu|u_{n}|^{2^*-2}u_{n}\phi- u_n\phi\log u^2_n\big)\mathrm{d}x=o_n(1)\|\phi\|_{W},\ \ \text{for each}\ \ \phi\in W.
\end{equation}
In particular, if we take $\phi=u_{n}$,
we then conclude $\{\lambda_{n}\}$ is bounded. Hence, up to a subsequence, we can assume $\lambda_{n}\to \lambda \in \R$.
Passing to the limit in \eqref{b411}, we get that  $u$ satisfies
\begin{equation*}
-\Delta u+\lambda u=u\log u^2+\mu|u|^{2^*-2}u,\ \ x\in \R^N.
\end{equation*}
It follows that $P(u)=0$, we conclude that
\begin{equation*}
\|\nabla u\|^2_2=\mu\|u\|^{2^*}_{2^*}+\frac{N}{2}c^2.
\end{equation*}
Set $v_n := u_n-u$, by the Br\'{e}zis-Lieb Lemma \cite{WM} and $P(u_n)=o_n(1)$, we deduce that
\begin{equation}\label{x3}
\int_{\R^N}|\nabla v_n|^2\mathrm{d}x=\mu\int_{\R^N}|v_n|^{2^*}\mathrm{d}x+o_n(1).
\end{equation}
Here we distinguish the following two cases
\begin{equation*}
\text{either}\quad \ (i) \ \|v_n\|^{2^*}_{2^*}\to 0\quad \text{or}\quad (ii)\ \|v_n\|^{2^*}_{2^*}\to l>0.
\end{equation*}
\vskip1mm
If $(ii)$ holds, we deduce from \eqref{x3} that
\begin{equation*}
\|\nabla v_n\|^2_2\ge \mu^{-\frac{N-2}{2}}\mathcal{S}^{\frac{N}{2}}+o_n(1).
\end{equation*}
Thus, we have
\begin{equation}\label{x4}
\begin{aligned}
m^-(c)+o_n(1)
&=\frac{1}{2}\|\nabla u_n\|^2_2-\frac{1}{2}\int_{\R^N}u^2_n\log u^2_n\mathrm{d}x
-\frac{\mu}{2^*}\int_{\R^N}|u_n|^{2^*}\mathrm{d}x+\frac{1}{2}c^2\\
&\ge \frac{1}{2}\|\nabla v_n\|^2_2-\frac{\mu}{2^*}\int_{\R^N}|v_n|^{2^*}\mathrm{d}x+I(u)\\
&\ge \frac{1}{N}\mu^{-\frac{N-2}{2}}\mathcal{S}^{\frac{N}{2}}+m^+(c),\\
\end{aligned}
\end{equation}
which contradicts with Lemma \ref{lem5.9}.
\vskip1mm

If $(i)$ holds, then $\{u_n\}\subset H^1_r(\R^N)$ converges strongly in $H^1_r(\R^N)$. It follows from \eqref{x4} that
$\int_{\R^N}u^2_n\log u^2_n\mathrm{d}x \to \int_{\R^N}u^2\log u^2\mathrm{d}x$, then we deduce $u_n\to u$ strongly in $W_r$.
By Lemma \ref{lem5.3}, we then conclude $u$ is a positive solution of \eqref{eq1.1}-\eqref{eq1.11}, which is radially symmetric and decreasing in $r=|x|$.
\qed

\vskip 0.08truein

{\bf Proof of Theorem \ref{th1.3}.}
The proof follows directly from Lemma \ref{Lem5.8} and Lemma \ref{lem5.10}.

\section{Proof of Theorem \ref{th1.4}}
In this section, we prove Theorem \ref{th1.4}.
In order to study the asymptotic behavior of ground state as $\mu\to 0$, we first introduce some facts for equations \eqref{eq1.1}-\eqref{eq1.11} with $\mu=0$.
Let $c>0$ be fixed, $(\lambda_0,w_0)\in \R \times W$ solves the following equation
\begin{equation}\label{v3}
\begin{cases}
-\Delta u+\lambda u=u\log u^2,\\
\int_{\R^N} u^2\mathrm{d}x=c^2.
\end{cases}
\end{equation}
Define
$$m_0(c):=\inf_{u\in S(c)}I_0(u),$$
where
\[
I_0(u)=\frac{1}{2}\int_{\R^N}|\nabla u|^2+u^2\mathrm{d}x-\frac{1}{2}\int_{\R^N}u^2\log u^2\mathrm{d}x.
\]
From \cite{DMS,CT}, by scaling, we obtain the unique solution of \eqref{v3} is
\begin{equation}\label{v4}
w_0(x):=c\pi^{-\frac{N}{4}}e^{-\frac{|x|^2}{2}},\quad \text{with} \quad \lambda_0:=-N-\log (c^{-2}\pi^{\frac{N}{2}}).
\end{equation}
Therefore $m_0(c)=\frac{c^2}{2}\big(N+1+\log (c^{-2}\pi^{\frac{N}{2}})\big)$ and $I(w_0)=m_0(c)$.

\begin{lem}\label{lem6.1}
Assume  $(\lambda_{\mu},u^+_{\mu})$ is the normalized ground state solution of \eqref{eq1.1}-\eqref{eq1.11} in Theorems \ref{th1.1}-\ref{th1.3},
then, up to a subsequence, we have
\begin{equation*}
u^+_{\mu}\to w_0\quad \text{strongly~in}\ \ W,
\end{equation*}
and $\lambda_{\mu}\to \lambda_0$ as $\mu \to 0^+$, where $(\lambda_0,w_0)$ is defined by \eqref{v4}. 
\end{lem}

\noindent{\bf Proof.}
Case $(i)$: if $(\lambda_{\mu},u^+_{\mu})$ is the normalized ground state solution of \eqref{eq1.1}-\eqref{eq1.11} obtained by Theorem \ref{th1.1},
we then claim that $\lim_{\mu\to 0}m(c)=m_0(c)$.
\vskip1mm
Indeed, if $\mu>0$, we have $\limsup_{\mu\to 0^+}m(c)\le m_0(c)$ since $I(w_0)<I_0(w_0)=m_0(c)$ for each $w_0\in S(c)$. Furthermore, for any $u\in S(c)$, one has
\begin{equation}\label{v61}
m_0(c)\le I_0(u)=I(u)+\frac{\mu}{p}\int_{\R^N}|u|^p\mathrm{d}x.\\
\end{equation}
Letting $\mu\to 0^+$, we then get $m_0(c)\le \liminf_{\mu\to 0^+} m(c)$. Thus, the claim hold.

\vskip1mm

Therefore, as $\mu\to 0^+$, we may assume that
\begin{equation}\label{v5}
I(u^+_{\mu})=\frac{1}{2}\int_{\R^N}|\nabla u^+_{\mu}|^2+|u^+_{\mu}|^2\mathrm{d}x-\frac{1}{2}\int_{\R^N}|u^+_{\mu}|^2\log|u^+_{\mu}|^2\mathrm{d}x
-\frac{\mu}{p}\int_{\R^N}|u^+_{\mu}|^p\mathrm{d}x\le m_0(c)+1.
\end{equation}
As in the proof of Lemma \ref{lem3.2}, \eqref{v5} gives that $u^+_{\mu}$ is bounded in $W$. Then, there exists $u_0\in W$ such that
$u^+_{\mu}\rightharpoonup u_0$ weakly in $W$, $u^+_{\mu} \to u_0$ strongly in $L^q(\R^N)$, $2\le q<2^*$, and $u^+_{\mu} \to u_0$ a.e. $x\in \R^N$. Thus,
\begin{equation*}
\begin{aligned}
m_0(c)&\le I_0(u_0)
\le \liminf_{\mu\to 0^+}\frac{1}{2}\int_{\R^N}\Big(|\nabla u^+_{\mu}|^2+|u^+_{\mu}|^2+A(|u^+_{\mu}|)\Big)\mathrm{d}x-\frac{1}{2}\int_{\R^N}B(|u_0|)\mathrm{d}x\\
&\le \liminf_{\mu\to 0^+}I(u^+_{\mu})+\frac{\mu}{p}\int_{\R^N}|u^+_{\mu}|^p\mathrm{d}x= m(c)+o(1).
\end{aligned}
\end{equation*}
Since $m(c)\to m_0(c)$, $I_0(u_0)=m_0(c)$, by the uniqueness of ground state solution of \eqref{v3}, we have $u_0=w_0$. Moreover,
$u^+_{\mu}\to w_0$ strongly in $W$, i.e. $u^+_{\mu}\to w_0$ strongly in $H^1(\R^N)$ and $\int_{\R^N}|u^+_{\mu}|^2\log|u^+_{\mu}|^2\mathrm{d}x
\to \int_{\R^N}|w_0|^2\log|w_0|^2\mathrm{d}x$ as $\mu \rightarrow 0^+$, and then $\lambda_{\mu}\to \lambda_0$ as $\mu \to 0^+$.
\vskip1mm
If $\mu<0$, we deduce that $m_0(c)=I_0(w_0)\le I(w_0)$ for $w_0\in S(c)$, and then $m_0(c)\le \liminf\limits_{\mu\to 0^-}m(c)$. For each $u\in S(c)$, one has
\begin{equation*}
\begin{aligned}
m(c)&\le I(u)=I_0(u)-\frac{\mu}{p}\int_{\R^N}|u|^p\mathrm{d}x.\\
\end{aligned}
\end{equation*}
Letting $\mu\to 0^-$, we get that $m(c)\le \liminf\limits_{\mu\to 0^-} m_0(c)$. This prove that $\lim\limits_{\mu\to 0^-}m(c)=m_0(c)$. In the same manner we obtain that $u^+_{\mu}\to w_0$ strongly in $W$ and $\lambda_{\mu}\to \lambda_0$ as $\mu \rightarrow 0^-$.

\vskip2mm
Cases $(ii)$, if $(\lambda_{\mu},u^+_{\mu})$ is the normalized ground state solution of \eqref{eq1.1}-\eqref{eq1.11} obtained by Theorems \ref{th1.2}-\ref{th1.3},
we claim that $\lim\limits_{\mu\to 0^+}m^+(c)=m_0(c)$.
By Lemma \ref{lem4.5} and Lemma \ref{lem5.7}, $m_c=m^+(c)$, we only need to prove $\lim\limits_{\mu\to 0^+}m_c=m_0(c)$.
From \eqref{z2} and $\eqref{eqs5.7-7}$, $m_c=\inf\limits_{u\in V_{k_0}}I(u)$, we then get $\limsup\limits_{\mu\to 0}m_c\le m_0(c)$ since $w_0\in V_{k_0}$.
Moreover, for $u\in V_{k_0}$, \eqref{v61} shows that $m_0(c)\le \liminf\limits_{\mu\to 0^+} m(c)$ as $\mu\to 0^+$.
Then,
\begin{equation*}
\begin{aligned}
I(u^+_{\mu})&= \Big(\frac{1}{2}-\frac{1}{p\gamma_p}\Big)\|\nabla u^+_{\mu}\|^2_2+\frac{1}{2}\int_{\R^N}A(|u^+_{\mu}|)-B(|u^+_{\mu}|)\mathrm{d}x-\frac{p}{2(p-2)}c^2,\\
\end{aligned}
\end{equation*}
we deduce $u^+_{\mu}$ is bounded in $W$. In the same manner, we can show that $u^+_{\mu}\to w_0$ strongly in $W$ and $\lambda_{\mu}\to \lambda_0$.
\qed


{\bf Proof of Theorem \ref{th1.4}.}
The proof directly follows from Lemma \ref{lem6.1}.

\section{The case $\alpha<0$}
In this section, we assume that $\alpha<0$ and $N\ge 2$. For convenience, we assume $\alpha=-1$.
Then the energy functional is defined by
\begin{equation*}
I(u)=\frac{1}{2}\int_{\R^N}|\nabla u|^2\mathrm{d}x+\frac{1}{2}\int_{\R^N}u^2\log u^2\mathrm{d}x
-\frac{\mu}{p}\int_{\R^N}|u|^{p}\mathrm{d}x-\frac{1}{2}\int_{\R^N}u^2\mathrm{d}x,\ \ u\in W_r,
\end{equation*}
and the Pohozaev identity is gived by
\begin{equation*}
P(u)=\|\nabla u\|^2_2+\frac{N}{2}c^2-\mu\gamma_p\|u\|^p_p=0.
\end{equation*}

{\bf Proof of Theorem \ref{th1.51}.}
We notice that if one of the following conditions holds
\begin{equation*}
(i)\ p:=\bar{p}=2+\frac{4}{N},\mu>0 \  \text{and} \ c<\Big(\frac{N+2}{N \mu}\Big)^{\frac{N}{4}}\big(C(N,\bar{p})\big)^{-\frac{N+2}{2}},
\quad (ii)\ \mu\le 0  \ \text{and} \ p>2;
\end{equation*}
then the fiber map $\Psi_{u}(t):=I(t\star u)$ is strictly increasing for each $u\in S(c)$, and thus $I$ does not have critical points on $S(c)$.
\vskip1mm
$(iii)$ By contradiction, we suppose that there exists a positive solution $u\in W_r$, by \cite[Radial
Lemma A.IV]{BL1}, $\lim\limits_{|x|\to \infty}u(x)=0$. Therefore, there exists $R_0>0$ large enough such that
\begin{equation*}
-\Delta u(x)=\big(-\lambda -\log u^2+\mu|u|^{p-2}\big)u(x)\ge u(x)\quad \text{for}~|x|>R_0.
\end{equation*}
By applying Theorem 2.1 and Theorem 2.8 of \cite{AB} with $f(s):=s$, we obtain that $-\Delta u\ge f(u)$
has no positive solution in any exterior domain of $\R^N$ if $\liminf_{s\to 0}s^{-\frac{N}{N-2}}f(s)>0$ for $N=3$ or
if $\lim_{s\to \infty}e^{as}f(s)=\infty$ for every $a>0$ and $N=2$.
\qed

\subsection{The existence of a global maximal on $\mathcal{P}_c$}
Define
\begin{equation*}
D:=\Big(\frac{Np\gamma_p}{2(2-p\gamma_p)}\Big)^{\frac{2-p\gamma_p}{2(p-2)}}\Big(\frac{\mu p\gamma^2_p}{2}C^p(N,p)\Big)^{-\frac{1}{p-2}}.
\end{equation*}
We next show that $\mathcal{P}_c$ is not empty if $c\geq D$. 
\begin{lem}\label{lem7.1}
Let $\mu>0$, $2<p<2+\frac{4}{N}$ and $N\ge 2$, 
we have
\begin{equation*}
\mathcal{P}_c\neq \emptyset \quad \text{if and only if}\quad c\ge D.
\end{equation*}
\end{lem}

\noindent{\bf Proof.}
For each $u\in S(c)$, define
\begin{equation*}
\phi_{u}(t):=P(t\star u)=e^{2t}\|\nabla u\|^2_2+\frac{N}{2}c^2-\mu\gamma_pe^{p\gamma_pt}\|u\|^p_p.
\end{equation*}
Therefore
\begin{equation*}
\inf_{t\in \R}\phi_{u}(t)=-\frac{2-p\gamma_p}{p\gamma_p}\Big[\frac{\mu p\gamma^2_p\|u\|^p_p}{2\|\nabla u\|^2_2}\Big]^{\frac{2}{2-p\gamma_p}}\|\nabla u\|^2_2
+\frac{N}{2}c^2.
\end{equation*}
By the Gagliardo-Nirenberg inequality, we obtain
\begin{equation*}
\inf_{t\in \R}\phi_{u}(t)\ge \frac{N}{2}c^2-\frac{2-p\gamma_p}{p\gamma_p}\Big(\frac{\mu p\gamma^2_p}{2}C^p(N,p)\Big)^{\frac{2}{2-p\gamma_p}}
c^{\frac{2p(1-\gamma_p)}{2-p\gamma_p}}.
\end{equation*}
If
\begin{equation*}
c<D:=\Big(\frac{Np\gamma_p}{2(2-p\gamma_p)}\Big)^{\frac{2-p\gamma_p}{2(p-2)}}\Big(\frac{\mu p\gamma^2_p}{2}C^p(N,p)\Big)^{-\frac{1}{p-2}},
\end{equation*}
we then get $\inf\limits_{u\in S(c)}P(u)>0$ and $\mathcal{P}_c=\emptyset$. Since the best constant in the Gagliardo-Nirenberg inequality is achieved, we assume that by $\bar{u}\in S(c)$, and then 
\begin{equation}\label{y2}
\inf_{u\in S(c)}P(u)=P(\bar{u})=\frac{N}{2}c^2-\frac{2-p\gamma_p}{p\gamma_p}\Big(\frac{\mu p\gamma^2_p}{2}C^p(N,p)\Big)^{\frac{2}{2-p\gamma_p}}
c^{\frac{2p(1-\gamma_p)}{2-p\gamma_p}}.
\end{equation}
If $c>D$, it follows from \eqref{y2} that $\inf\limits_{u\in S(c)}P(u)<0$. Since $\lim\limits_{t\to +\infty}\phi_{u}(t)=+\infty$, by continuity of $\phi_{u}$,
there exists $\bar{t}\in \R$ such that $\phi_{u}(\bar{t})=0$, i.e. $\mathcal{P}_c\neq\emptyset$. Moreover, if $c=D$, then $P(\bar{u})=0$, and $\mathcal{P}_c\neq\emptyset$ holds.
\qed

\begin{lem}\label{lem7.2}
Let $2<p<2+\frac{4}{N}$ and $N\ge 2$, if $c\ge D$, we have
\begin{equation*}
\inf_{u\in \mathcal{P}_c}I(u)=-\infty.
\end{equation*}
\end{lem}
\noindent{\bf Proof.}
From Lemma \ref{lem7.1}, if $c\ge D$, then $\mathcal{P}_c\neq\emptyset$. For any $u\in S(c)$ and $P(u)\le 0$, since $\lim\limits_{t\to +\infty}P(t\star u)=+\infty$,
by the geometry of $I(t\star u)$, we deduce there exists $t\ge 0$ such that $P(t\star u)=0$ and $I(t\star u)\le I(u)$.
Therefore, we only need to prove that there exists a sequence $\{u_n\}\subset S(c)$ with $P(u_n)\le 0$ and $I(u_n)\to -\infty$ as $n\to \infty$.
We choose $\eta>0$ and take $u\in \mathcal{C}^{\infty}_0(\R^N)$, $u\ge 0$ with $\|u\|^2_2=c^2-\frac{\eta}{2}$ and $P(u)< 0$.
We also choose a $v\in \mathcal{C}^{\infty}_0(\R^N)$, $v\ge 0$ with $\|v\|^2_2=\frac{\eta}{2}$. We define
\begin{equation*}
u_n(x):=u(x)+\Big(\frac{1}{n}\Big)^{\frac{N}{2}}v\Big(\frac{1}{n}(x-n^2e_1)\Big)=u(x)+v_n(x),
\end{equation*}
where $e_1=(1,0,\cdots,0)$ and $n>0$ is chosen sufficiently large so that the supports of $u$ and $v_n$ are disjoints. We then obtain
\begin{equation*}
\begin{aligned}
P(u_n)&=\int_{\R^N}|\nabla(u+v_n)|^2\mathrm{d}x+\frac{N}{2}c^2-\mu \gamma_p\int_{\R^N}|u+v_n|^p\mathrm{d}x\\
&=\int_{\R^N}|\nabla u|^2\mathrm{d}x+\frac{N}{2}c^2-\mu \gamma_p\int_{\R^N}|u|^p\mathrm{d}x+o_n(1)\le 0,\\
\end{aligned}
\end{equation*}
since $\|\nabla v_n\|^2_2\to 0$ and $\|v_n\|^{p}_{p}\to 0$ as $n\to \infty$. Taking into account that $\int_{\R^N}u^2_n\log u^2_n\mathrm{d}x\to -\infty$, we then get
\begin{equation*}
I(u_n)\to-\infty,\quad \text{as}~n\to \infty.
\end{equation*}
Thus, the lemma is proved.
\qed

\vskip 0.12truein

Define
$$k_0:=\frac{p\gamma_pN}{2(2-p\gamma_p)}c^2.$$

\begin{lem}\label{lem7.3}
Assume that $2<p<2+\frac{4}{N}$ and $N\ge 2$. If $P(u)\le 0$ and $\|\nabla u\|^2_2=k_0$, then
\begin{equation*}
c\ge D=\left(\frac{1}{\mu\gamma_p C^p(N,p)}\Big(\frac{N}{2-p\gamma_p}\Big)^{\frac{2-p\gamma_p}{2}}\big(\frac{p\gamma_p}{2}\big)^{-\frac{p\gamma_p}{2}}\right)^{\frac{1}{p-2}}.
\end{equation*}
Moreover,  If $P(u)<0$ and $\|\nabla u\|^2_2=k_0$, then $c>D$.
\end{lem}

\noindent{\bf Proof.}
Since $P(u)\le 0$, by Gagliardo-Nirenberg inequality, we have
\begin{equation*}
\|\nabla u\|^2_2+\frac{N}{2}c^2\le \mu\gamma_p C^{p}(N,p)c^{p(1-\gamma_p)}\|\nabla u\|^{p\gamma_p}_2.
\end{equation*}
The above equality holds only for the best constant in the Gagliardo-Nirenberg inequality is achieved.
Let $\|\nabla u\|^2_2=\frac{p\gamma_pN}{2(2-p\gamma_p)}c^2$, we get
\begin{equation*}
c\ge D=\left(\frac{1}{\mu\gamma_p C^p(N,p)}\Big(\frac{N}{2-p\gamma_p}\Big)^{\frac{2-p\gamma_p}{2}}\big(\frac{p\gamma_p}{2}\big)^{-\frac{p\gamma_p}{2}}\right)^{\frac{1}{p-2}}.
\end{equation*}
 The lemma is completed.
\qed

\begin{lem}\label{lem7.4}
Let $2<p<2+\frac{4}{N}$ and $N\ge 2$, if $c= D$, then
\begin{equation*}
M(c):=\sup_{u\in \mathcal{P}_c\cap H^1_r(\R^N)}I(u)
\end{equation*}
is achieved. 
\end{lem}
\noindent{\bf Proof.}
From Lemma \ref{lem7.1}, if $c\ge D$, then $\mathcal{P}_c\neq\emptyset$. Let $u\in \mathcal{P}_c$, by Gagliardo-Nirenberg inequality, we have
\begin{equation}\label{y3}
\|\nabla u\|^2_2+\frac{N}{2}c^2\le \mu\gamma_p C^{p}(N,p)c^{p(1-\gamma_p)}\|\nabla u\|^{p\gamma_p}_2,
\end{equation}
then $\|\nabla u\|^2_2$ is bounded from above. 
Moreover, we deduce that
\begin{equation}\label{y4}
\begin{aligned}
I(u)&= \Big(\frac{1}{2}-\frac{1}{p\gamma_p}\Big)\|\nabla u\|^2_2+\frac{1}{2}\int_{\R^N}B(|u|)-A(|u|)\mathrm{d}x-\frac{p}{2(p-2)}c^2\\
&\le \Big(\frac{1}{2}-\frac{1}{p\gamma_p}\Big)\|\nabla u\|^2_2+K_qC^{q}(N,q)c^{q(1-\gamma_q)}\|\nabla u\|^{q\gamma_q}_2-\frac{p}{2(p-2)}c^2,
\end{aligned}
\end{equation}
where $q\gamma_q\le 2$. Therefore, $I$ restricted to $\mathcal{P}_c$ is bounded from above, and $M(c)<+\infty$.
Let $\{u_n\}\subset \mathcal{P}_c\cap H^1_r(\R^N)$ be a maximizing sequence for $I$ on $M(c)$, combining \eqref{y3} and \eqref{y4},
$u_n$ is bounded in $H^1_r(\R^N)$ and $A(|u_n|)$ is bounded in $L^1(\R^N)$. It follows from Lemma \ref{lem2.1} that $u_n$ is bounded in $W_r$.
From Lemma \ref{lem2.2}, there exists $u\in S(c)\cap H^1_r(\R^N)$
 such that $u_n\rightharpoonup u$ weakly in $W_r$ and $u_n\to u$ strongly in $L^s(\R^N)$, $2\le s<2^*$. Moreover, $P(u)\le 0$.
By a direct calculation, if $P(u)=0$, then $I(u)\le M(c)$. Define
\begin{equation*}
F(u)=I(u)-\frac{1}{p\gamma_p}P(u)=\Big(\frac{1}{2}-\frac{1}{p\gamma_p}\Big)\|\nabla u\|^2_2
+\frac{1}{2}\int_{\R^N}u^2\log u^2\mathrm{d}x-\frac{p}{2(p-2)}c^2.
\end{equation*}
Since $A$ is lower semicontinuous for the weak convergence in $W_r$, we have $I(u)\ge M(c)$.
By contradiction, $P(u)<0$. Taking into account that
\begin{equation*}
F(t\star u)=\Big(\frac{1}{2}-\frac{1}{p\gamma_p}\Big)e^{2t}\|\nabla u\|^2_2+\frac{tN}{2}\int_{\R^N}u^2\mathrm{d}x
+\frac{1}{2}\int_{\R^N}u^2\log u^2\mathrm{d}x-\frac{p}{2(p-2)}c^2,
\end{equation*}
if $\frac{d}{dt}|_{t=0}F(t\star u)=0$, we have $\|\nabla u\|^2_2=k_0$, which contradicts with Lemma \ref{lem5.3}.

\vskip1.2mm

We now prove $\frac{d}{dt}|_{t=0}F(t\star u)=0$ by distinguishing two cases.
\vskip1mm
{\bf Case 1.} We first assume that $\frac{d}{dt}|_{t=0}F(t\star u)<0$.
Since $\lim\limits_{t\to-\infty}P(t\star u)=\frac{N}{2}c^2>0$ and $P(u)<0$,
there exists $t_0\in\R$ such that $t_0<0$ and $t_0\star u\in \mathcal{P}_c$. It follows that $P(t\star u)< 0$ for $t\in (t_0,0]$.
We claim that $\frac{d}{dt}F(t\star u)<0$ for all $t\in (t_0,0]$. In fact, by Lemma \ref{lem7.3}, $\frac{d}{dt}F(t\star u)\neq0$ for $\bar{t}\in (t_0,0]$.
If $\frac{d}{dt}|_{\bar{t}}F(t\star u)>0$ for $\bar{t}\in (t_0,0]$, then by continuity, there exists $t_1\in [\bar{t},0]$ such that
$\frac{d}{dt}|_{t_1}F(t\star u)=0$, which contradicts with Lemma \ref{lem7.3}. Therefore, we have
\begin{equation*}
 F(t_0\star u)=I(t_0\star u)>F(u)=I(u)-\frac{1}{p\gamma_p}P(u)\ge M(c),
\end{equation*}
this is a contradiction with the definition of $M(c)$.
\vskip1mm
{\bf Case 2.} Next, we assume that $\frac{d}{dt}|_{t=0}F(t\star u)>0$. Since $\lim\limits_{t\to +\infty}P(t\star u)=+\infty$, there exists
$t_2>0$ such that $P(t_2\star u)=0$ and $P(t\star u)<0$ for all $t\in [0,t_2)$. Thus, similar to the proof of Case 1,
we have $\frac{d}{dt}|_{\bar{t}}F(t\star u)>0$ for $t\in [0,t_2)$. We also have $F(t_2\star u)>M(c)$, which leads to the same contradiction.
\vskip1mm
In conclusion, we deduce that $P(u)=0$. Therefore, we have $I(u)=M(c)$.
\qed

\begin{lem}\label{lem7.5}
Let $p\in(2,2+\frac{8}{N(N+2)})\bigcup (2+\frac{8}{N(N+2)},2+\frac{4}{N})$, if $c= D$, for any maximizer of $M(c)$ is a critical point of $I$ restricted to $S(c)$.
\end{lem}

\noindent{\bf Proof.}
We prove that $\mathcal{P}_{c}\cap H^1_r(\R^N)$ is a smooth manifold of codimension 1 in $S(c)$.
By Lemma \ref{lem7.4}, we deduce there exists $u\in \mathcal{P}_{c}\cap H^1_r(\R^N)$ such that
$$I(u)=\max\limits_{\mathcal{P}_{c}\cap H^1_r(\R^N)}I.$$
Then there exist $\lambda_1,\lambda_2\in \R$ such that $dI(u)=\lambda_1 u+\lambda_2 d P(u)$, and
\begin{equation*}
-(1+2\lambda_2)\Delta u +\lambda_1u+u\log u^2=\left(1+p\gamma_p\lambda_2\right)\mu|u|^{p-2}u.\\
\end{equation*}
The Pohozaev identity for the above equation is
\begin{equation}\label{y5}
(1+2\lambda_2)\|\nabla u\|^2_2+\frac{N}{2}c^2=\left(1+p\gamma_p\lambda_2\right)\mu\gamma_p\|u\|^p_p.
\end{equation}
Therefore, by combining \eqref{y5} with $P(u)=0$, we have
\begin{equation}\label{x5}
2\lambda_2\|\nabla u\|^2_2=p\gamma^2_p\lambda_2\mu\|u\|^p_p.
\end{equation}
If $\lambda_2\neq0$, it follows from $P(u)=0$ that
\begin{equation}\label{x7}
\|\nabla u\|^2_2+\frac{N}{2}c^2\le \mu\gamma_p C^{p}(N,p)c^{p(1-\gamma_p)}\|\nabla u\|^{p\gamma_p}_2.
\end{equation}
The equality in \eqref{x7} holds only for the best constant in the Gagliardo-Nirenberg inequality is reached, from \cite{WeM}, 
$u$ satisfies
\begin{equation}\label{x6}
\|\nabla u\|^2_2=\|u\|^2_2.
\end{equation}
Combining \eqref{x5} with \eqref{x6}, we have $p=2+\frac{8}{N(N+2)}$, which is a contradiction.
Therefore, the strict inequality in \eqref{x7} holds,
which contradicts with Lemma \ref{lem7.3}, then $\lambda_2=0$.

In conclusion, $\mathcal{P}_{c}$ is a smooth manifold of codimension 1 on $S(c)$.
\qed

\vskip 0.12truein

{\bf Proof of Theorem \ref{th1.5}.}
By Lemma \ref{lem7.4}, we know that there exists $u\in W_r$ such that $I(u)=\sup_{\mathcal{P}_c\cap H^1_r(\R^N)}I$. From Lemma \ref{lem7.5}, $u$ is a critical point of $I$ restricted to $S(c)$. Moreover, it directly follows from Theorem \ref{th1.51} $(iii)$ that $u$ is non-positive.

\vskip 0.12truein

{\bf Data availability statement: }
Data sharing is not applicable to this article as no new data were created or
analyzed in this study.


\end{document}